\documentclass[11pt, a4paper,leqno]{amsart}
\usepackage{amsmath,amsthm,amscd,amssymb,amsfonts, amsbsy}
\usepackage{latexsym}
\usepackage{exscale}
\usepackage{yhmath}
\usepackage{mathrsfs}
\usepackage{enumitem}
\usepackage[utf8]{inputenc}
\usepackage{datetime}

\usepackage[colorlinks,citecolor=red,hypertexnames=false, breaklinks]{hyperref}

\calclayout
\allowdisplaybreaks

\newdateformat{monthdayyeardate}{%
	\monthname[\THEMONTH]~\THEDAY, \THEYEAR}
\newcommand{\X}{\mathbf{X}}
\renewcommand{\emptyset}{\mbox{\textup{\O}}}

\newcommand{\normH}[1]{{\vert\kern-0.25ex\vert\kern-0.25ex\vert #1 
		\vert\kern-0.25ex\vert\kern-0.25ex\vert}}

\newcommand{\normHB}[1]{{\Big\vert\kern-0.25ex\Big\vert\kern-0.25ex\Big\vert #1 
		\Big\vert\kern-0.25ex\Big\vert\kern-0.25ex\Big\vert}}

\def\Xint#1{\mathchoice
	{\XXint\displaystyle\textstyle{#1}}%
	{\XXint\textstyle\scriptstyle{#1}}%
	{\XXint\scriptstyle\scriptscriptstyle{#1}}%
	{\XXint\scriptscriptstyle\scriptscriptstyle{#1}}%
	\!\int}
\def\XXint#1#2#3{{\setbox0=\hbox{$#1{#2#3}{\int}$}
		\vcenter{\hbox{$#2#3$}}\kern-.5\wd0}}
\def\fint{\Xint-}

\def\Xint#1{\mathchoice
	{\XXint\displaystyle\textstyle{#1}}%
	{\XXint\textstyle\scriptstyle{#1}}%
	{\XXint\scriptstyle\scriptscriptstyle{#1}}%
	{\XXint\scriptscriptstyle\scriptscriptstyle{#1}}%
	\!\int}
\def\XXint#1#2#3{{\setbox0=\hbox{$#1{#2#3}{\int}$}
		\vcenter{\hbox{$#2#3$}}\kern-.5\wd0}}

\newtheorem*{step1}{Step 1}

\newcommand{\dv}{\mbox{div}\hspace{0.1cm}}
\def \R{ \mathbb{R}}

\theoremstyle{plain}
\newtheorem{theorem}{Theorem}[section]
\newtheorem{lemma}{Lemma}[section]

\newtheorem{proposition}{Proposition}[section]

\newtheorem{claim}{Claim}[section]

\theoremstyle{definition}
\newtheorem{definition}{Definition}[section]

\theoremstyle{remark}
\newtheorem{remark}{Remark}[section]

\DeclareMathOperator*{\esssup}{ess\,sup}

\newcommand{\Lm}{\mathscr{L}}
\newcommand{\E}{\mathbf{X}\times \mathbf{Z}}
\newcommand{\rn}{\mathbb{R}^{n-1}}
\newcommand{\rnp}{\mathbb{R}^{n}_+}

\newcommand{\va}{\varLambda}

\newcommand{\rnn}{\mathbb{R}^{n}}

\newcommand{\kk}{\mathrm{k}}
\newcommand{\g}{\mathrm{g}}

\newcommand{\Z}{\mathbf{Z}}
\newcommand{\ds}{\star\star}

\numberwithin{equation}{section}

\begin{document}
	\allowdisplaybreaks
	\author{Gael Y. Diebou}
	\address{Gael Y. Diebou\\
		University of Toronto, Department of Mathematics
		\\
40 St George Street, Toronto, ON
		\\
M5S 2E4, Canada}
	\email{gaely.diebou@utoronto.ca}
	
	\title[]{Existence and regularity of steady-state solutions of the Navier-Stokes equations arising\\ from irregular data}
	\date{\today}
	
	\thanks{The author was partially funded by the Fields Institute for Research in Mathematical Sciences}

	
	\begin{abstract}
		We analyze the forced incompressible stationary Navier-Stokes flow in  $\rnp$, $n> 2$. Existence of a unique solution satisfying a global integrability property measured in a scale of tent spaces is established for small data in homogenous Sobolev space with $s=-\frac{1}{2}$ degree of smoothness. The velocity field is shown to be locally H\"{o}lder continuous while the pressure belongs to $L^p_{loc}$ for any $p\in (1,\infty)$. Our approach is based on the analysis of the inhomogeneous Stokes system for which we derive a new solvability result involving Dirichlet data in Triebel-Lizorkin classes with negative amount of smoothness and is of independent interest. 
	\end{abstract}
	
	\maketitle

	
	\bigskip
	\section{Introduction }\label{s:1}
	The steady state (forced) incompressible Navier-Stokes equations in a domain $\Omega\subset \rnn$, $n\geq 2$ is the following system 
	\begin{align}\label{eq:NS-eq}\tag{\textbf{NS}}
		\begin{cases}
			-\Delta u+\nabla \pi+u\cdot\nabla u =F\hspace{0.2cm}\mbox{in}\hspace{0.2cm}\Omega\\
			\dv u=0\hspace{0.2cm}\mbox{in}\hspace{0.2cm}\Omega
		\end{cases}
	\end{align}
	where $u:\Omega\rightarrow \rnn$ is the unknown velocity field, $\pi:\Omega\rightarrow \R$ is the unknown scalar pressure and $F:\Omega\rightarrow \rnn$ is a given external force. This system is supplemented by the boundary condition 
	\begin{equation}\label{bc}
		u=f\hspace{0.2cm}\mbox{on}\hspace{0.2cm}\partial\Omega
	\end{equation} 
	where $f=(f_1,...,f_n)$ is a prescribed vector field satisfying (in the case $\Omega$ smooth bounded) the compatibility condition $\displaystyle\int_{\partial\Omega}f\cdot \textbf{N}d\sigma=0$
	with $\textbf{N}=(N_1,\cdots,N_{n})$ being the outer unit normal vector at the boundary.
	
	Probably, the first striking result regarding the solvability of the Dirichlet problem for the Navier-Stokes equations was obtained by Leray \cite{Le}. In a smooth bounded three dimensional domain, he showed the existence of a weak solution $(u,\pi)\in W^{1,2}(\Omega)\times L^2(\Omega)$ provided $f\in W^{1/2,2}(\partial\Omega)$ and $F\in W^{-1,2}(\Omega)$. Existence of generalized weak solutions of \eqref{eq:NS-eq}-\eqref{bc} in $\Omega\subset \mathbb{R}^3$, those are $(u,\pi)\in W^{1,q}(\Omega)\times L^q(\Omega)$,  $q\in (2,\infty)$  is a consequence of the important work by Cattabriga \cite{Catta}. Without some extra condition on the data, the regularity of weak solutions to \eqref{eq:NS-eq}-\eqref{bc} is still far from being fully understood. Available contributions are mostly dimension dependent. Indeed, in the physically  most relevant dimensions $n=2,3$; any weak solution is smooth (see e.g. \cite{Se} in the case of nonhomogeneous smooth data). In four dimensions,  the $L^p$-regularity theory of weak solutions with zero Dirichlet data was established in  \cite{Ger}. When $\Omega=\mathbb{R}^5$ or the torus $\mathbb{T}^5$, it is known due to Struwe \cite{Str} that \eqref{eq:NS-eq} admits a solution $(u,\pi)\in W^{2,p}_{loc}\times W^{1,p}_{loc}$ for any $p<\infty$ whenever the external force is smooth and compactly supported. This result was partially extended to the case $\Omega=\mathbb{T}^n$ and very recently to $\Omega=\mathbb{R}^n$ with $n\in [5,15]$ in \cite{Fr-Ru} and \cite{YYY}, respectively. As for uniqueness of weak solutions, it seems that a smallness assumption on given data is necessary and a recent result by Luo \cite{Luo} predicts that this condition cannot be dropped.
	
	There have been growing interest in recent years in the analysis of the Navier-Stokes equations subject to low regularity data. By this, we mean that $f$ satisfies a weaker regularity  than that needed in order for generalized weak solutions to exist. In such a context, does problem \eqref{eq:NS-eq}-\eqref{bc} admit a solution? In the affirmative case, what are the qualitative properties of such solutions?
	Prescribing  boundary value with low regularity forces one to consider a notion of solution weaker than weak solution. A good candidate, roughly speaking is obtained by testing \eqref{eq:NS-eq} against a suitable  divergence-free smooth vector field and performing two successive integration by parts taking \eqref{bc} into account. This idea to the best of the author's knowledge first appeared in \cite{Am} and the resulting formulation gives rise to the so-called very weak solution. When $\Omega\subset \rnn$, $n=2,3$ is a $C^2$ regular bounded domain, the author in \cite{M-P} constructed such a solution in $L^{2n/(n-1)}(\Omega)$ provided $f\in L^2(\partial\Omega)$ (with arbitrarily large norm) and $F\in W^{-1,2n/(n-1)}(\Omega)$. Existence,  uniqueness and regularity of very weak solutions $(u,\pi)$ in the class $L^q(\Omega)\times W^{-1,q}(\Omega)$ have been obtained in \cite{FGS,GSS} under certain smallness conditions on  $f\in W^{-1/q,q}(\partial\Omega)$ and $F\in W^{-1,r}(\Omega)$ with $1<r\leq q<\infty$, $1/r\leq 1/n+1/q$. These results were generalized in \cite{Kim} wherein the author gave a complete solvability theory for very weak solutions of \eqref{eq:NS-eq}-\eqref{bc}. In particular, refining the definition of very weak solutions and using some ideas from the preceding references, the author showed the existence of a very weak solution $(u,\pi)$ in $L^{n}(\Omega)\times W^{-1,n}(\Omega)$ for arbitrary large data $f\in W^{-1/n,n}(\partial\Omega)$ and forcing term $F$ for $n=3,4$. In two dimensions, he proved the existence of a solution $(u,\pi)\in L^{q_0}(\Omega)\times W^{-1,q_0}(\Omega)$, $2<q_0<3$. Moreover, he investigated the regularity of these solutions and also derived uniqueness results under suitable smallness assumptions. The existence theory for very weak solutions in unbounded domains (the half-space, exterior domains, etc.) seems to be more subtle. In general, similar methods as those employed for instance in \cite{Kim} which rely on duality arguments and functional analytic tools cannot be carried out. We refer the reader to \cite{Finn} for an interesting discussion pertaining to generalized (weak) solutions -- we point out however, some recent existence results for the linear Stokes problem in half-space domain \cite{Fsa} and in exterior domains \cite{KKP}.  
	
	This paper aims at establishing the solvability theory for \eqref{eq:NS-eq}-\eqref{bc} in the case $\Omega=\rnp$ (with possible  adaptation of the ideas to the case of bounded smooth domains) by means of novel ideas. The techniques employed here are new and complement those introduced in \cite{DK} for the analysis of elliptic "critical" problems subject to low regularity data.
	In addition, they can be invoked to study similar questions for other semilinear elliptic problems. Assuming $\Omega=\rnp$, we seek for a velocity field of the form $u=v+w$ where $v$ solves the linear Stokes equation with Dirichlet data $f$ while $w$ solves the inhomogeneous Stokes problem with zero boundary data and source term $F-u\cdot\nabla u$. Odqvist \cite{Odqvist} proved  that $v$ assumes an integral representation, it is the Stokes extension of $f$ to $\rnp$ (see Section \ref{s:2}). We then look for the data $f$ in a large class of distributions on $\rn$ for which $v$ is well-defined and has $f$ as  trace in a suitable sense.  On one hand,
	\eqref{eq:NS-eq} is scaling (and translation) invariant with respect to the maps
	\begin{align*}
		u_{\lambda}(x)=\lambda u(\lambda x),\quad \pi_{\lambda}(x)=\lambda^2 \pi(\lambda x),\quad \lambda>0
	\end{align*} 
	for appropriately rescaled external force and boundary data. On the other hand, we want to have $u$ in the (local) Lebesgue space $L^2_{loc}$ in order to make sense of the equation. From these observations, we are led to the consideration of $v$ in $T^{2(n-1),2}_{s_2}$, $s_2=-\frac{1}{2(n-1)}$ a scale of weighted tent spaces introduced in \cite{CMS,HMM}. Thanks to the continuous characterization of homogeneous Trieblel-Lizorkin spaces \cite{Trie}, we expect that $v$ must have a distributional trace $f$ in the homogeneous negative Sobolev space $\dot{H}^{-1/2,2(n-1)}(\rn)$.  By the same token, the pressure $\pi$ is sought for in the weighted tent space $T^{2(n-1),2}_{\sigma_2}$, $\sigma_2=-\frac{3}{2(n-1)}$ (see below for the definition of weighted tent spaces). These expectations are actually confirmed using some special properties of the Stokes kernel (smoothness, cancellation), see Proposition \ref{prop:Stokes-Triebel}. Tent spaces naturally arise in the analysis of linear elliptic equations and systems, see e.g. \cite{AAlex,HMM} and references therein. We quote the recent work \cite{DK} where these spaces are used in the context of nonlinear systems. 
	
	Our main result states that there exists a unique solution of \eqref{eq:NS-eq}-\eqref{bc} in a suitable framework under a smallness condition on $f\in\dot{H}^{-1/2,2(n-1)}(\rn)$ and the external force. A more general statement involving Dirichlet data in homogeneous Triebel-Lizorkin spaces is obtained. In both cases, the constructed solution satisfies a global integrability property, expressed in terms of a weighted tent norm and it is further shown that this solution enjoys a better regularity locally. More precisely, the velocity field $u$ belongs to $C^{0,\alpha}_{loc}$ and the pressure $\pi$ is an element of $L^p_{loc}$ for any $p\in (1,\infty)$. These features are derived from the pointwise decay rate of the velocity field near the boundary. In order to derive all these results,  we study the inhomogeneous Stokes problem in $\rnp$ (which plays a fundamental role in the analysis of \eqref{eq:NS-eq}-\eqref{bc} when the flow takes place in an exterior region, a channel or a pipe) and establish key estimates of the solution for prescribed data in the homogeneous Triebel-Lizorkin class with negative amount of smoothness.   
	\subsection{Tent spaces and functional settings}
	Throughout, a point $x\in \rnp$ will typically be denoted by $(x',x_n)$, $x'\in \rn$ and $x_n>0$. For $R>0$, $B_R(x')$ is the  closed ball in $\rn$ with radius $R>0$ and center at $x'\in \rn$. Given $\alpha>0$, define the cone (nontangential region)  with vertex at $x'\in \rn$ by  \[\varGamma_{\alpha}(x'):=\{(y',y_n)\in \rnp:|x'-y'|<\alpha y_n\}.\]
	When $\alpha=1$,  $\varGamma_1(x')$ will instead be denoted as  $\varGamma(x')$.
	Given a ball $B=B_R(x')$, we denote by  $T(B)=B_R(x')\times (0,2R)$ the Carleson box over $B_R(x')$. 
	For $q\in [1,\infty)$, consider the functionals $\mathcal{A}_q$, $\mathscr{C}_q$ defined for $F$ measurable in $\rnp$ by
	\begin{equation}\label{conical-funct}
		\mathcal{A}_q F(x')=\bigg(\iint_{\varGamma(x')}|F(y',y_n)|^qy_n^{-n}dy'dy_n\bigg)^{1/q},\hspace{0.12cm} \mathcal{A}_{\infty} F(x')=\esssup_{(y',y_n)\in \varGamma(x')} |F(y',y_n)| 
	\end{equation}
	\begin{equation}\label{Carl-funct}
		\mathscr{C}_qF(x')=\sup_{x'\in B}\bigg(\dfrac{1}{|B|}\iint_{T(B)}|F(y',y_n)|^qy_n^{-1}dy'dy_n\bigg)^{1/q}
	\end{equation}
	where the supremum is taken over all balls containing $x'$.	The membership of each of these functionals in Lebesgue spaces gives rise to a scale of function spaces  first introduced by Coifman, Meyer and Stein \cite{CMS}.   The tent space $T^{p,q}$ with $p,q\in [1,\infty)$ collects all functions $F\in L^q_{loc}(\rnp)$  for which $\mathcal{A}_qF\in L^p(\rn)$. We equip this space with the norm   
	\begin{equation}
		\|F\|_{T^{p,q}}:=\|\mathcal{A}_qF\|_{L^p(\rn)}.
	\end{equation}
	When $p=\infty$, the space $T^{\infty,q}$ is defined using the Carleson functional $\mathscr{C}_q$ by
	\[T^{\infty,q}=\{F\in L^q_{loc}(\rnp):\mathscr{C}_qF\in L^{\infty}(\rn)\}.\]
	The space $T^{\infty,q}$ is intrinsically linked to Carleson measures. In fact, it can alternatively be interpreted as the space of functions $F\in L^{q}_{loc}(\rnp)$ for which $d\mu(y',y_n)=|F|^qy_n^{-1}dy'dy_n$ is a Carleson measure in $\rnp$. Such measures have some interesting application in the study of nonlinear PDEs with critical growth in the gradient, see \cite{DK}. Note that at the endpoint case $p=q=\infty$, $T^{p,q}$ may be identified to $L^{\infty}(\rnp)$. 
	For $s \in \mathbb{R}$, we say that $F:\rnp\rightarrow \R$ belongs to the  weighted tent space \cite{Alex,HMM},  which we denote by $T^{p,q}_{s}$ if 
	\[(y',y_n)\mapsto y_n^{-(n-1)s}F(y',y_n)\in T^{p,q}.\] 
	It can  easily be verified that $\|F\|_{T^{p,q}_{s}}:=\big\|y_n^{(1-n)s}F\big\|_{T^{p,q}}$  defines a norm on $T^{p,q}_{s}$ and is a Banach space. In particular, the space $T^{p,q}_{s_q}$, $s_q=-\frac{1}{q(n-1)}$  contains the space of functions in $L^q(\rnp)$ with compact support as a dense subspace. This property, together with the completeness of $T^{p,q}_{s_q}$  follows from Lemma \ref{lem:equiv-L^q} below.
	\begin{lemma}\label{lem:equiv-L^q} Let $K$ be a compact set in $\rnp$ and assume that $F\in T^{p,q}_{s_q}$ for $p,q\in [1,\infty)$. Then
		\begin{align}\label{equiv-L^q}
			C_1\|\mathbf{1}_KF\|_{T^{p,q}_{s_q}}\leq \|F\|_{L^q(K)}\leq C_2\|F\|_{T^{p,q}_{s_q}}  
		\end{align}
		where the constant $C_1,C_2$ only depend on $p,q,n$ and $K$.
	\end{lemma}
	Out of convenience, we defer the proof of Lemma \ref{lem:equiv-L^q} to the Appendix. 
	\begin{remark}\label{rmk:ind-aperture}
		The change of aperture of the cone does not modify the tent norm in the sense that if one defines the $\alpha$-conical functional $\mathcal{A}^{\alpha}_q$ by 
		\begin{align*}
			\mathcal{A}^{\alpha}_q F(x'):=\bigg(\iint_{\varGamma_{\alpha}(x')}|F(y',y_n)|^q(\alpha y_n)^{-(n-1)}\frac{dy'dy_n}{y_n}\bigg)^{1/q},\quad \alpha>0    
		\end{align*}
		then 
		\begin{align}
			\label{Aq}\|\mathcal{A}^{\alpha}_qF\|_{L^{p}(\rn)}\approx \|\mathcal{A}^{\beta}_qF\|_{L^p(\rn)}
		\end{align}
		where the implicit constant depends on $p,q$ and $\alpha,\beta\in (0,\infty)$. See \cite[Proposition 4, p. 309]{CMS} which remains valid for $q\neq 2$. 
	\end{remark}
	Moreover, for $s_1,s_2\in \mathbb{R}$ such that $s_1>s_2$ and $1\leq p_1<p_2\leq \infty$, $q\in (0,\infty)$ the following continuous embedding holds (see \cite[Lemma 2.19]{Alex})
	\begin{equation}\label{embedding-tents}
		T^{p_1,q}_{s_1}\subset T^{p_2,q}_{s_2}    
	\end{equation}
	provided $s_2-s_1=\dfrac{1}{p_2}-\dfrac{1}{p_1}$. Recall H\"{o}lder's inequality in weighted tent spaces. 
	\begin{lemma}\label{lem:Holder}
		Let $p_i,q_i\in [1, \infty]$ and $s_i\in \R$,   $i\in \{0,1,2\}$ such that $\displaystyle\sum_{i=1}^21/p_i=1/p_0$ and $\displaystyle\sum_{i=1}^21/q_i=1/q_0$  with the convention $1/\infty=0$. If $f\in T^{p_1,q_1}_{s_1}$ and $g\in T^{p_2,q_2}_{s_2}$, then $fg\in T^{p_0,q_0}_{s_0}$ and it holds that
		\begin{equation}\label{Holder-tent}
			\|fg\|_{T^{p_0,q_0}_{s_0}}\leq C \|f\|_{T^{p_1,q_1}_{s_1}} \|g\|_{T^{p_2,q_2}_{s_2}} 
		\end{equation}
		provided $s_0=s_1+s_2$.
	\end{lemma}
	This lemma can be proved via factorization of tent spaces \cite[Theorem 3.4]{Hu}.  It is long-established that there is an intrinsic connection between weighted tent spaces and Triebel-Lizorkin spaces which we now recall their definitions. 
	
	Let us denote by $\mathcal{S}(\rn)$ the class of Schwartz (smooth rapidly decreasing) functions on $\rn$  and $\mathcal{S}'(\rn)$ its topological dual space endowed with the weak-$\star$ topology. Define the space \[\dot{\mathcal{S}}(\rn)=\{f\in \mathcal{S}(\rn)\hspace{0.1cm}\big|\hspace{0.1cm}\int z^{\gamma}f(z)dz=0,\hspace{0.1cm}\forall \,\gamma\in \mathbb{N}^{n} \}\]
	which inherits the topology of $\mathcal{S}(\rn)$ as subspace. This space may be identified with the space of Schwartz functions whose Fourier transforms vanish together with all their derivatives at the origin. Its dual space is denoted by $\dot{\mathcal{S}}'(\rn)$.
	Let $\varphi \in \mathcal{S}(\rn)$ such that
	\begin{equation*}
		\varphi(\xi)=
		\begin{cases}
			1\hspace{0.2cm}\mbox{if}\hspace{0.2cm}|\xi|\leq 1\\
			0\hspace{0.2cm}\mbox{if}\hspace{0.2cm}|\xi|> 2.
		\end{cases}    
	\end{equation*}
	Let $\psi(\xi)=\varphi(\xi)-\varphi(2\xi)$ so that  
	\begin{equation*}
		\sum_{j\in \mathbb{Z}}\psi(2^{-j}\xi)=1,\hspace{0.2cm}\forall \,\xi\in \mathbb{R}^{n-1}\setminus\{0\}.   
	\end{equation*}
	Let us denote by $\mathcal{F}f$ the Fourier transform of $f$ on $\rn$ and by $\dot{\Delta}_j=\mathcal{F}^{-1}(\psi(2^{-j}\cdot)\mathcal{F}\cdot)$ the homogeneous Littlewood-Paley operator. For
	$s\in \mathbb{R}$ and $p,q\in [1,\infty)$, we say that a distribution $f\in \dot{\mathcal{S}}'(\rn)$ belongs to the Triebel-Lizorkin space $\dot{F}^{s}_{p,q}(\rn)$ if 
	\begin{equation*}
		\|f\|_{\dot{F}^{s}_{p,q}(\rn)}=\bigg\|\bigg[\sum^{\infty}_{j=-\infty}\big(2^{js}\big|\dot{\Delta}_jf\big|\big)^q\bigg]^{1/q}\bigg\|_{L^p(\rn)}<\infty.    
	\end{equation*}
	This space is complete and equivalent to the homogeneous Sobolev space
	$\dot{H}^{s,p}(\rn)$
	whenever $q=2$ and $1<p<\infty$. Moreover, for $1\leq q_1,q_2\leq \infty$ and $s_1,s_2\in \R$ with $s_2<s_1$ we have the continuous inclusion
	\begin{equation*}
		\dot{F}^{s_1}_{p_1,q_1}(\rn)\subset \dot{F}^{s_2}_{p_2,q_2}(\rn)     
	\end{equation*}
	provided $p_1,p_2\in (1,\infty)$ with $s_1-s_2=\frac{n-1}{p_1}-\frac{n-1}{p_2}>0$.
	\begin{remark}\label{rmk:disccharac}
		Let $\phi\in \mathcal{S}(\rn)$ with Fourier transform $\widehat{\phi}$ supported in $\{1\leq |\xi|\leq 2\}$ and $\sum_{k\in \mathbb{Z}}\big(\widehat{\phi}(2^{-k}\xi)\big)^2=1$ for all $\xi\in \rn\setminus \{0\}$. For each $k\in \mathbb{Z}$, let $\phi_k=2^{k(n-1)}\phi(2^k\cdot)$ be the dyadic dilation of $\phi$. Whenever $\phi_k\ast f$,  $f\in \mathcal{S}'(\rn)$ is a function we define the discrete  Peetre maximal function
		\begin{equation}\label{disPee}
			\phi^{\star}_kf(x')=   \sup_{y'\in \rn}|(\phi_k\ast f)(x'-y')|(1+2^{k}|y'|)^{-\lambda},\quad x'\in \rn.
		\end{equation}
		This object can also be used to characterize function spaces of Triebel-Lizorkin and Besov type. In particular, if $\lambda>\max\{(n-1)/p,(n-1)/q\}$, then 
		\begin{equation}\label{discharac}
			\bigg\|\bigg(\sum_{k\in \mathbb{Z}}(2^{ks}\phi_k^{\star}f)^{q}\bigg)^{1/q}\bigg\|_{L^p(\rn)} \leq C\|f\|_{\dot{F}^s_{p,q}}  
		\end{equation}
		where the constant $C$ is independent of $f$. It is interesting that in \eqref{discharac}, the kernel $\phi$ needs not be smooth, see for instance \cite{BC}. 
		
	\end{remark}
	\begin{definition}\label{defn:space-1}
		For $n> 2$, $p\in (1,\infty)$ and $q\in [1,p)$  we define $\X_{p,q}$ as the space of vector fields $u:\rnp\rightarrow \rnn$ satisfying $\|u\|_{\X_{p,q}}<\infty$ and $\Z_{p,q}:=\{\pi:\rnp\rightarrow \mathbb{R}\hspace{0.1cm}{\big |}\hspace{0.1cm}\|\pi\|_{\Z_{p,q}}<\infty\}$ 
		where
		\begin{align*}
			\|u\|_{\X_{p,q}}&=\sup_{x_n>0}x_n^{\frac{1}{q}+\frac{n-1}{p}}\|u(\cdot,x_n)\|_{L^{\infty}(\rn)}+\|u\|_{T_{s_q}^{p,q}},\quad s_q=-\frac{1}{q(n-1)}
		\end{align*}
		and 
		\begin{align*}
			\|\pi\|_{\Z_{p,q}}:= \|\pi\|_{T_{\sigma_{q}}^{p,q}},\hspace{0.2cm} \sigma_q=-\frac{1+1/q}{n-1}.
		\end{align*}	
	\end{definition}
	Note that either $\|\cdot\|_{\X_{p,q}}$ or $\|\cdot\|_{\Z_{p,q}}$ defines a norm on $\X_{p,q}$ and $\Z_{p,q}$ respectively, they are Banach spaces. For convenience, when $p=2(n-1)$ and $q=2$, we will specially denote the spaces $\X_{p,q}$ and $\Z_{p,q}$ by $\X$ and $\Z$, respectively.  
	
	\subsection{Main results}
	Our first result deals with the well-posedness theory. In what follows, the dimension is assumed larger or equal to $3$ unless otherwise stated.
	
	\begin{theorem}\label{thm:unforced}Assume that $F=0$. 
		Then the Dirichlet problem \eqref{eq:NS-eq}-\eqref{bc} has a solution $(u,\pi)$ in $\X\times \hspace{0.1cm}\Z$ which is unique in a small closed ball provided the data $f$ has a sufficiently small $[\dot{H}^{-\frac{1}{2},2(n-1)}(\rn)]^{n}$-norm.
	\end{theorem}
	
	In presence of the forcing term, our main finding reads as follows.
	
	\begin{theorem}\label{thm1}Let $\eta\in (1,2)$ and $\eta<\tau<2(n-1)$ verify the condition $\dfrac{1}{\eta}+\dfrac{n-1}{\tau}=3$. 
		There exists $\varepsilon>0$ such  that for every $f\in [\dot{H}^{-\frac{1}{2},2(n-1)}(\rn)]^{n}$ and $F\in \X_{\tau,\eta}$ satisfying $\|f\|_{\dot{H}^{-\frac{1}{2},2(n-1)}(\rn)}+\|F\|_{\X_{\tau,\eta}}< \varepsilon$, Eq. \eqref{eq:NS-eq}-\eqref{bc} admits a solution $(u,\pi)$ in $\X\times \hspace{0.1cm}\Z$ which is the only one among those satisfying the condition
		$\|u\|_{\X}+\|\pi\|_{\Z}\leq 2\varepsilon$. 
	\end{theorem}
	
	Existence of solutions in $\X_{p,q}\times \Z_{p,q}$ for any $2< q<\infty$ is a consequence of Theorem \ref{thm1} together with an improved regularity result. In more details, the statement reads as follows. 
	
	\begin{theorem}	\label{thm:solvability-persistance} 
		Let $\eta,\tau$ and $\varepsilon>0$ be as in Theorem \ref{thm1}, take $2<q<p<\infty$, $\eta_1\in (1,q)$ and $\eta_1<\tau_1<p$.  Given $f=(f_1,...,f_n)$ in $\big[\dot{H}^{-\frac{1}{2},2(n-1)}\cap  \dot{F}^{s}_{p,q}(\rn)\big]^{n}$ and $F\in \X_{\tau,\eta}\cap \X_{\tau_1,\eta_1}$, there exists $\varepsilon^{\star}\in (0,\varepsilon)$ such that if $\|f\|_{\dot{H}^{-\frac{1}{2},2(n-1)}(\rn)}+\|F\|_{\X_{\tau,\eta}}< \varepsilon^{\star}$ then there exists a solution $(u,\pi)$ of \eqref{eq:NS-eq}-\eqref{bc} in the space $\X_{p,q}\times \Z_{p,q}$ which is unique in the ball
		\begin{align*}
			B_{2\varepsilon^{\star}}(\mathbf{0},0)=\{(u,\pi)\in \X_{p,q}\times \Z_{p,q}: \|u\|_{\X_{p,q}}+\|\pi\|_{\Z_{p,q}}\leq 2\varepsilon^{\star}\}
		\end{align*}	
		provided $\frac{1}{\eta_1}+\frac{n-1}{\tau_1}=2+\frac{1}{q}+\frac{n-1}{p}$ and $s=-\frac{1}{q}$.
	\end{theorem}
	The uniqueness of the pressure as claimed in the previous results should be understood up to an additive constant.  Let us now record the following regularity result which arises as a consequence of the local boundedness property of the velocity field and the elliptic regularity theory for the Stokes system \cite{Galdi}.
	\begin{theorem}\label{thm2}
		If $(u,\pi)\in \X\times \Z$ is the solution of problem \eqref{eq:NS-eq}-\eqref{bc} constructed in Theorem \ref{thm:unforced} or \ref{thm1}, then $(u,\pi)\in [C^{0,\alpha}_{loc}(\rnp)]^{n}\times L^{p}_{loc}(\rnp)$ for some $\alpha\in (0,1)$ and every $p\in (1,\infty)$. 
	\end{theorem}	
The regularity $s=-1/2$ is special in the analysis of \eqref{eq:NS-eq}-\eqref{bc}. In fact, the space $T^{2(n-1),2}_{s_2}$ (resp. the boundary class $\dot{H}^{-1/2,2(n-2)}(\rn)$) is one scale  within the one parameter  family of spaces $\bigg\{T_{\frac{n-1-p}{(n-1)p}}^{p,2}\bigg\}_{2<p<\infty}$ (resp. the space $\dot{F}^{s}_{p,2}(\rn)$ with $s=\frac{n-1}{p}-1<0$) which are all admissible according to the scaling symmetry. However, if $p>2$ is fixed and $H=u\otimes u\in T^{p/2,1}_{\frac{2(n-1-p)}{(n-1)p}}$, then in order to estimate the solution in $\X_{p,2}$ (the weighted $L^{\infty}$-norm to be more precise) by the corresponding norm of $H$, the requirement $p\leq 2(n-1)$, is necessary, at least from our arguments (see for instance the proof of Proposition \ref{prop:Green-pot-bound}) so that $s=-1/2$ is the endpoint. This is a particular case of the natural condition $\omega(y_n)=y_n^{-2s}\in \mathcal{A}_2$, $\mathcal{A}_2$ being the Muckenhoupt weight class.
	
	\begin{remark}
		It is important to emphasize that the solution $(u,\pi)$ found in the preceding theorems is in fact a weak solution of \eqref{eq:NS-eq}-\eqref{bc} $($compare with the notion of regular solution defined in \cite{YYY}$)$  in the sense that 
		\[u\in W^{1,p}_{loc}\quad \mbox{and}\quad \pi\in L^p_{loc},\] 
		for any $p<\infty$.		
		Pertaining to Theorem \ref{thm:unforced}, the construction of a solution $(u,\pi)$ in the framework $L^2(\mathbb{R}^2_+)\times L^2(\mathbb{R}^2_+,x_ndx)$, that is when $n=2$, seems challenging -- the main difficulty in the analysis comes from the critical regularity of the nonlinear term $u\cdot \nabla u$. 
	\end{remark}
	
	\section{Auxiliary results}\label{s:2}
	This section is devoted to the analysis of the Dirichlet problem for the following system
	\begin{align}\label{eq:S-eq}\tag{\textbf{S}}
		\begin{cases}
			-\Delta u+\nabla \pi=F-\dv H \hspace{0.2cm}\mbox{in}\hspace{0.2cm}\rnp\\
			\dv u=0\hspace{0.2cm}\mbox{in}\hspace{0.2cm}\rnp\\
			u=f\hspace{0.2cm}\mbox{on}\hspace{0.2cm}\partial\rnp
		\end{cases}
	\end{align}
	with given vector fields $f,F$ and tensor $H$. Our goal is to prove that \eqref{eq:S-eq} admits a  solution $(u,\pi)$ in the target space $\X_{p,q}\times \Z_{p,q}$ whose norm can be estimated by the norms of $f$, $F$ and $H$ in suitable function spaces. To this end, for better readability we simply separate the study into two parts: the homogeneous case ($f=0$) and the inhomogeneous case ($F=0$, $H=\textbf{0}$). 
	
	\subsection{Homogeneous Stokes system and linear estimates}
	Consider the Stokes operator $L_{S}$ acting on pair of functions $(u,\pi)\in [\mathscr{D}'(\rnn)]^{n}\times \mathscr{D}'(\rnn)$, $n\geq 2$ and given by 
	\[L_{S}(u,\pi)=\big(-\Delta u_1+\partial_1\pi,\cdots,-\Delta u_{n}+\partial_{n}\pi, \sum_{i=1}^{n}\partial_iu_i\big).\]
	Let us denote by $\mathcal{M}_{n\times n}(\mathcal{S}'(\rnn))$ the collection of all $n\times n$ matrices with coefficients in $\mathcal{S}'(\rnn)$. A fundamental solution of the Stokes operator $L_{S}$ in $\rnn$ is a pair $(\mathbf{E},\mathbf{b})$ with $\mathbf{E}=(E_{ij})^{n}_{i,j=1}\in \mathcal{M}_{n\times n}(\mathcal{S}'(\rnn))$ and $\textbf{b}=(b_1,...,b_{n})\in [\mathcal{S}'(\rnn)]^{n}$ satisfying coordinate-wise the equations
	\begin{align*}
		\begin{cases}
			-\Delta E_{ij}+\partial_i b_j=\delta_{ij}\delta\hspace{0.2cm}\mbox{in}\hspace{0.2cm}\mathcal{S}'(\rnn),\hspace{0.2cm}\forall\,\,i,j\in \{1,...,n\}\\
			\displaystyle\sum_{k=1}^{n}\partial_kE_{kj}=0\hspace{0.2cm}\mbox{in}\hspace{0.2cm}\mathcal{S}'(\rnn),\hspace{0.2cm} \forall\,\,j\in \{1,...,n\}.
		\end{cases}
	\end{align*} 
	For $n\geq 3$, one may apply the Fourier transform to both sides of each of the above equations after which one deduces the  explicit expressions
	\begin{equation}\label{eq:fund-sols}
		E_{ij}(x)=\frac{1}{2\omega_{n-1}}\bigg[\frac{1}{(n-2)}\frac{\delta_{ij}}{|x|^{n-2}}+\frac{x_ix_j}{|x|^{n}}\bigg], \hspace{0.1cm} b_j=\frac{1}{\omega_{n-1}}\frac{x_j}{|x|^{n}};\hspace{0.12cm} i,j\in \{1,...,n\}
	\end{equation}  
	defined for $x\in \rnn\setminus\{0\}$ where $\omega_{n-1}$ is the surface measure of the unit sphere in $\rnn$ centered at zero. Details of computations leading to \eqref{eq:fund-sols} can be found in \cite[Chap. 10]{DM}. On the other hand, when $n=2$, $\textbf{E}$ and $\textbf{b}$ assume the following forms
	\begin{equation}
		E_{ij}(x)=\frac{1}{4\pi}\bigg[\frac{x_ix_j}{|x|^{2}}-\delta_{ij}\log |x|\bigg], \hspace{0.1cm} b_j=\frac{1}{2\pi}\frac{x_j}{|x|^{2}};\hspace{0.12cm} x\in \mathbb{R}^2\setminus\{0\},\hspace{0.1cm} i,j\in \{1,2\}.
	\end{equation}
	Now, let us consider the homogeneous Stokes system 
	\begin{align}\label{eq:Hom-Stokes}
		\begin{cases}
			-\Delta u+\nabla \pi=0\hspace{0.2cm}\mbox{in}\hspace{0.2cm}\rnp\\
			\dv u=0\hspace{0.2cm}\mbox{in}\hspace{0.2cm}\rnp\\
			u=f\hspace{0.2cm}\mbox{on}\hspace{0.2cm}\partial\rnp.
		\end{cases}
	\end{align}
	With the convolution being understood in a component wise sense, define
	\begin{equation}\label{eq:linear-int-eqs}\mathcal{H}f(x',x_n)=(\mathcal{K}_{x_n}\ast f)(x'),\quad\mathcal{E}f(x',x_n)=(\kk_{x_n}\ast f)(x')
	\end{equation}
	where $$\mathcal{K}_{x_n}(x')=(K_{ij}(x',x_n))_{1\leq i,j\leq n}\hspace{0.2cm} \mbox{and} \hspace{0.2cm}\kk_{x_n}(x')=(\kk_1(x',x_n),...,\kk_{n}(x',x_n))$$ are commonly referred to as the Odqvist kernels \cite{Odqvist} -- each entry of the tensors assuming an explicit form in terms of $\mathbf{E}$ and $\mathbf{b}$ via the formulas
	\begin{equation}\label{eq:sol-Homo-Stokes}K_{ij}(x)=2(\partial_{n}E_{ij}+\partial_jE_{in}-\delta_{jn}b_i)=-\frac{2n}{\omega_{n-1}}\frac{x_nx_ix_j}{|x|^{n+2}} 
	\end{equation}
	and 
	\begin{equation}\label{eq:sol-press}
		\kk_j(x)=4\partial_jb_{n}=\frac{4}{\omega_{n-1}}\partial_j\frac{x_n}{|x|^{n}}.
	\end{equation}
	For the full derivation of these kernels, the interested reader may consult the articles \cite{Odqvist, Solo}. Note that if $f$ belongs to the weighted Lebesgue space $L^1\big(\rn,\frac{dx'}{(1+|x'|)^{n}}\big)$, then $u=\mathcal{H}f$ and $\pi=\mathcal{E}f$ are both meaningful as absolutely convergent integrals and $(u,\pi)$ is the unique solution of Eq. \eqref{eq:Hom-Stokes} decaying at infinity. The Stokes extension of a tempered distribution may not be meaningful in general (as another tempered distribution). Harmonic Poisson extensions of Schwartz distributions have been studied by H. Triebel  \cite{Trie} -- they characterize almost all scale of Triebel-Lizorkin spaces on $\rn$  via tent spaces. In particular, the following equivalence holds true   
	\begin{equation}\label{equiv}
		\bigg\|\mathcal{A}_q[x_n^{m-s}\partial^m_{x_n}\mathcal{P}_{x_n}f]\bigg\|_{L^{p}(\rn)}\sim\|f\|_{\dot{F}^{s}_{p,q}(\rn)},\quad 1\leq p,q<\infty
	\end{equation}
	where $m$ is a large nonnegative integer, $m>s$ and $\mathcal{P}_{x_n}f$ denotes the convolution of $f$ with the kernel $\mathcal{P}_{x_n}(x')=\frac{2}{\omega_{n-1}}x_n(|x'|^{2}+x_n^{2})^{-\frac{n}{2}}$.  It is not straightforward to see that the convolution in \eqref{equiv} is meaningful given $f\in \dot{F}^{s}_{p,q}(\rn)$ since the considered kernel is non-smooth.  Stein \cite{Stein} introduced the notion of bounded distributions (i.e. $f\ast \phi$ $\in L^{\infty}(\rn)$ for any $\phi\in \mathcal{S}(\rn)$) and showed that if $f$ is a bounded distribution then $\mathcal{P}_{x_n}f$ for $x_n>0$ fixed, is a bounded smooth function in $\rnp$. This notion was recently extended in \cite{BC} to distributions of finite growth.
	\begin{definition}
		A tempered distribution $f$ in $\rn$ is said to be of finite growth $\gamma\geq 0$ if $f\ast \phi=O(|x'|^{\gamma})$ as $|x'|\rightarrow \infty$ whenever $\phi \in \mathcal{S}(\rn)$.    
	\end{definition}
	By this, we then understand bounded distributions as those with growth $\gamma=0$. Given a kernel $\psi$ such that $(1+|\cdot|)^{\gamma}\psi\in L^1(\rn)$, its convolution with a tempered distribution of finite growth makes sense as a tempered distribution. This idea was used in \cite{BC} to establish a non-smooth characterization of homogeneous Triebel-Lizorkin and Besov spaces via two essential arguments. Any $f\in \dot{F}^s_{p,q}(\rn)$ is up to an additive polynomial a distribution of growth $\gamma>s-(n-1)/p$ and $f\ast \psi$, $\psi\in \mathcal{A}_{\varLambda,m,r}$ is  a bounded continuous function using the Calder\'{o}n reproducing formula for elements in $\dot{F}^s_{p,q}(\rn)$, see \cite[Theorem 3.1]{BC}. Here, for $\varLambda\geq 0$, $m,r\in \R$, $\mathcal{A}_{\varLambda,m,r}$ is the class of functions in $L^1(\rn)$ such that $\widehat{\psi}\in C^{n+[\varLambda]}(\rn\setminus\{0\})$ and satisfies for all multi-indices $\alpha$ with $|\alpha|\leq n+[\varLambda]$ 
	\begin{equation}\label{cancelation-smootheness}
		|\partial^{\alpha}\widehat{\psi}(\xi')|=O(|\xi'|^{r-|\alpha|})\hspace{0.2cm}\mbox{as}\hspace{0.2cm} |\xi'|\rightarrow 0 \hspace{0.2cm}\mbox{and}\hspace{0.2cm} |\partial^{\alpha}\widehat{\psi}(\xi')|=O(|\xi'|^{-(n-1)-m})\hspace{0.2cm}\mbox{as}\hspace{0.2cm} |\xi'|\rightarrow \infty.
	\end{equation}
	From these arguments, the question regarding the meaning of the convolution in \eqref{equiv} is then settled  as each of the kernels involved satisfies \eqref{cancelation-smootheness}. 
	
	Now, we would like to investigate the validity of a one sided-estimate in \eqref{equiv} but having the Poisson kernel replaced by those in \eqref{eq:sol-Homo-Stokes} and \eqref{eq:sol-press}, respectively.
	Since the results of this paper deal with spaces of negative smoothness indexes and for simplicity sake, let us  assume from now on that $s$ is negative.  The following lemma collects some useful properties of $K_{jk}$ and $\kk_j$.
	\begin{lemma}\label{lem:Fourierkernels}
		Set $\sigma_n=\frac{2n}{\omega_{n-1}}$ and consider the kernels 
		\begin{equation*}
			P_0(x')=\frac{-\sigma_n}{(|x'|^2+1)^{\frac{n+2}{2}}},\,\,\, P_j(x')=  \frac{-\sigma_nx_j}{(|x'|^2+1)^{\frac{n+2}{2}}},\,\,\, Q_{jk}=\dfrac{-\sigma_n x_jx_k}{(|x'|^2+1)^{\frac{n+2}{2}}},\,\,\, j,k=1,...,n-1.  
		\end{equation*}
		Then $P_0$, $P_j$ and $Q_{jk}$ belong to $\bigcap\limits_{1\leq p\leq \infty}L^p(\rn)$. Moreover, if $\mathcal{F}_{x'}$ denotes the Fourier transform in $\rn$, then
		\begin{equation}
			(\mathcal{F}_{x'}P_0)(\xi')=|\xi'|e^{-|\xi'|}-e^{-|\xi'|};\,\,\, 
			(\mathcal{F}_{x'}P_j)(\xi')=i\xi_j'e^{-|\xi'|},\quad \xi'\in \rn.   
		\end{equation}
		\begin{equation}
			(\mathcal{F}_{x'}Q_{jk})(\xi')=-\delta_{jk}e^{-|\xi'|}+\frac{\xi'_{j}\xi'_k}{|\xi'|}e^{-|\xi'|},\,\,\,\xi'\in \rn\setminus\{0\}.   
		\end{equation}
	\end{lemma}	
	We know according to \cite[Corollary 2.3]{BC} that after possibly subtracting a suitable polynomial to $f\in \dot{F}^{s}_{p,q}(\rn)$, the resulting expression is a bounded distribution. Moreover, either of the kernels $P_0$, $P_j$ and $Q_{jk}$ belongs to $L^1(\rn)$ and satisfies the cancellation and smoothness condition \eqref{cancelation-smootheness}. As a consequence and by application of \cite[Theorem 3.1]{BC}, the convolutions in \eqref{eq:linear-int-eqs} are not only well-defined distributions but also bounded continuous functions. 
	\begin{remark}
		There is another simple argument showing that the distribution in \eqref{eq:linear-int-eqs} are in fact bounded functions whenever $f\in \dot{F}^{s}_{p,q}(\rn)$. Assume that $K\in L^1(\rn)$ is such that $\widehat{K}$ is rapidly decreasing and smooth except at the origin. Take $\eta\in \mathcal{S}(\rn)$ with $\widehat{\eta}=1$ near  the origin, one can write $K_{x_n}f=\eta_{x_n}\ast f \ast K_{x_n} + \chi_{x_n}\ast f$  where $\chi\in \mathcal{S}(\rn)$ is such that $\widehat{\chi}=(1-\widehat{\eta})\widehat{K}$ and conclude that $K_{x_n}f$ is bounded since  $f$ is a bounded distribution. 
	\end{remark}
	Another application of Lemma \ref{lem:Fourierkernels} is the following
	\begin{proposition}\label{prop:Stokes-Triebel}
		Let $s\in \R$, $s<0$ and $p,q\in (1,\infty)$. Take $\varLambda=\max\{(n-1)/q,(n-1)/p\}$ and let $r,m\in \R$ with $r>s$, $m+s>\lambda$ for some $\lambda\in (\varLambda,[\varLambda]+1)$. Assume further that $\psi\in \mathcal{A}_{\varLambda,m,r}$ and set $\psi_{x_n}(x')=x_n^{-(n-1)}\psi(x'/x_n)$, $(x',x_n)\in\rnp$. If $f\in \dot{F}^s_{p,q}(\rn)$, then there exists $C>0$ independent of $f$ such that
		\begin{equation}\label{eq:nonsmooth-est}
			\|\psi_{x_n}f\|_{T_{s/(n-1)}^{p,q}}\leq C\|f\|_{\dot{F}^s_{p,q}(\rn)},\quad \psi_{x_n}f=\psi_{x_n}\ast f.    
		\end{equation}
	\end{proposition}
	The proof of this result as well as that of Lemma \ref{lem:Fourierkernels} are both postponed to the Appendix. The previous proposition is interesting in that it provides a continuous characterization of Triebel-Lizorkin spaces by non-smooth kernels. In fact, $\psi$ can be either of the kernels in Lemma \ref{lem:Fourierkernels} - each of them belongs to the class $\mathcal{A}_{\varLambda,m,r}$ for suitable choices of parameters. Another special choice of $\psi$ is the function  $\frac{2}{\omega_{n-1}}(1+|\cdot|^2)^{-n/2}$ in $\rn$, which  corresponds to a particular case of the harmonic characterization in \eqref{equiv}. However, we do not claim the reverse inequality in \eqref{eq:nonsmooth-est} under the prescribed assumptions on $\psi$. 
	
	As discussed above, $\mathcal{H}f$ and $\mathcal{E}f$ for any $f\in \dot{F}^{-1/q}_{p,q}(\rn)$ are bounded continuous functions in the tangential variable. Our next lemma shows more, it establishes their precise decay rates as well as their membership to some weighted tent spaces.
	\begin{lemma}\label{lem:lin-est} Let $1<q<p<\infty$. There exists a constant $C:=C(n,p,q)>0$ such that 
		\begin{align}\label{eq:lin-est}
			\|\mathcal{H}f\|_{\X_{p,q}}+\|\mathcal{E}f\|_{\Z_{p,q}}+\sup_{x_n>0}x_n^{1+1/q+(n-1)/p}\|\mathcal{E}f(\cdot,x_n)\|_{L^{\infty}(\rn)}\leq C\|f\|_{\dot{F}^{-1/q}_{p,q}(\rn)}	
		\end{align}
		for all $f\in [\dot{F}^{-1/q}_{p,q}(\rn)]^{n}$ where $\mathcal{H}f$ and $\mathcal{E}f$ are defined as in \eqref{eq:linear-int-eqs}.	
	\end{lemma}
	We state two more auxiliary results which will be useful in the demonstration of Lemma \ref{lem:lin-est}. 
	\begin{lemma}[Averaging Lemma]\label{lem:Aver}
		Assume that $F\in L^q(\rnp)$, $q\geq 1$. We have 
		\begin{align*}
			\int_{\rn}\iint_{\varGamma(x')}|F(y',y_n)|^q\frac{dy'dy_n}{y_n^{n-1}}dx'=c\int_{\rnp}|F(y)|^qdy  
		\end{align*}
		where $c>0$ only depends on $n$, the dimension.
	\end{lemma}
	The proof of this identity follows from a simple application of Fubini–Tonelli's Theorem. In what follows, for a subset $S\subset \rnp$, we set $E(S)=\{x'\in \rn: S\cap \varGamma(x')\neq \emptyset\}$. 
	\begin{lemma}\label{lem:E_S} Let $K\subset \rnp$ be a compact set. Then $E(K)$ is open, its Lebesgue measure $|E(K)|$ is finite. Moreover, there exists a constant $\theta>0$ depending on $K$ only such that $|E(K)|\leq C\theta^{n-1}$ for some constant $C:=C(n)>0$.
	\end{lemma}
	\begin{proof}
		Let $x'\in E(K)$, there exists $(y',y_n)\in \rnp$ with $ (y',y_n) \in K$ and  $y'\in B_{y_n}(x')$. Putting $R=y_n-|x'-y'|>0$, it plainly follows that $B(x',R)\subset E(K)$. Moving on, we remark that $E(K)$ is actually bounded. Fix $b>1$ and observe that since $K$ is compact, one may assume without loss of generality that $K=Q\times [\ell(Q),b\ell(Q)]$ where $Q$ is a closed cube in $\rn$ with side-length  $\ell(Q)$. Hence, $|E(K)|\leq |B(x',b\ell(Q))|\leq C(b\ell(Q))^{n-1}$, $x'\in Q$. 
	\end{proof}
	Now we are ready to prove Lemma \ref{lem:lin-est}.
	\begin{proof}[Proof of Lemma \ref{lem:lin-est}]
		By a direct computation, we have for all $(x',x_n)\in \rnp$, $$\kk_j(x',x_n)=\begin{cases}
			-\frac{4n}{\omega_{n-1}}\frac{x_nx_j}{|x|^{n+2}}=2x^{-n}_nP_j(x'/x_n)\hspace{0.2cm}\mbox{if}\hspace{0.2cm}j=1,...,n-1\\
			2\partial_{x_n}\mathcal{P}_{x_n}(x')\hspace{3.2cm}\mbox{if}\,\,\,\,\,\,j=n
		\end{cases}$$ 
		where $P_j$ is as in Lemma \ref{lem:Fourierkernels}. Now set $(\overline{u},\overline{\pi})=(\mathcal{H}(f),\mathcal{E}(f))$ and for $x\in \rnp$ fixed, use the interior estimate for the linear Stokes problem (see e.g. \cite{Sim}) together with Lemma \ref{lem:Aver} to write
		\begin{align}
			\nonumber			|\overline{\pi}(x)|^q&\leq c|B_{x_n/2}(x)|^{-1}\int_{B_{x_n/2}(x)}|\overline{\pi}(y',y_n)|^qdy'dy_n\\
			\nonumber		&\leq c|B_{x_n/2}(x)|^{-1}\sum_{j=1}^{n}\int_{B_{x_n/2}(x)}|\kk_j(\cdot,y_n)\ast f_j|^qdy'dy_n\\
			\nonumber		&\leq c|B_{x_n/2}(x)|^{-1}\bigg(\int_{B_{x_n}(x')\times [x_n/3,2x_n]}|\partial_{y_n}\mathcal{P}_{y_n}f_n|^qdy'dy_n+\\
			\nonumber   & \hspace{4cm} \sum_{j=1}^{n-1}\int_{B_{x_n}(x')\times [x_n/3,2x_n]}|y_n^{-n}P_{j}(\cdot/y_n)\ast f_j|^qdy'dy_n\bigg)\\
			\nonumber	&\leq cx_n^{-n-q}\bigg(\int_{E(B_{x_n}(x')\times[x_n/3,2x_n])}\iint_{\varGamma(z')}|y_n^{1+1/q}\partial_{y_n}\mathcal{P}_{y_n}f|^q\frac{dy'dy_n}{y_n^{n}}dz'+\\
			\nonumber &\hspace{3cm} \sum_{j=1}^{n-1}\int_{E(B_{x_n}(x')\times[x_n/3,2x_n])}\iint_{\varGamma(z')}|y_n^{1/q}(P_{j})_{y_n}f_j|^q\frac{dy'dy_n}{y_n^{n}}dz'\bigg)\\	
			\nonumber &\leq cx_n^{-(1+\frac{1}{q}+\frac{n-1}{p})q}\bigg(\big\|\mathcal{A}_q[y_n^{1+\frac{1}{q}}\partial_{y_n}\mathcal{P}_{y_n}f_n]\big\|^{q}_{L^{p}(\rn)}+\sum_{j=1}^{n-1}\big\|\mathcal{A}_q[y_n^{1/q}(P_{j})_{y_n}f_j]\big\|^{q}_{L^{p}(\rn)}\bigg)\\
			\label{eq:ptw-pressure}	&\leq cx_n^{-(1+\frac{1}{q}+\frac{n-1}{p})q}\|f\|^{q}_{\dot{F}^{-1/q}_{p,q}(\rn)}.
		\end{align} 
		Observe that we have used H\"{o}lder's inequality and Lemma \ref{lem:E_S} to get the penultimate estimate and \eqref{equiv} and Proposition \ref{prop:Stokes-Triebel} for the last.  From the above remark on the kernel $\kk_j$, the estimate on $\mathcal{E}f$ in $T^{p,q}_{\sigma_q}$ is a consequence of the extrinsic characterization \eqref{equiv} (applied with $m=1$) and Proposition \ref{prop:Stokes-Triebel}. The same observation pertaining to the velocity field gives $\|\overline{u}\|_{T^{p,q}_{s_q}}\leq C\|f\|_{\dot{F}^{-1/q}_{p,q}(\rn)}$. To see this, we simply write
		\begin{equation}\label{eq:Kjk-fcts}
			K_{jk}(\cdot,x_n)=x_n^{-(n-1)}\begin{cases}
				Q_{jk}(x_n^{-1}\cdot)\hspace{0.2cm}\mbox{if}\hspace{0.2cm} j,k=1,...,n-1\\
				P_{k}(x_n^{-1}\cdot)\hspace{0.4cm}\mbox{if}\hspace{0.2cm} j=n,\,k=1,...,n-1\\
				P_0(x_n^{-1}\cdot)\hspace{0.4cm}\mbox{if}\hspace{0.2cm} j=k=n
			\end{cases}    
		\end{equation}
		and apply Proposition \ref{prop:Stokes-Triebel}. 
		It then remains to establish the bound
		\begin{equation}\label{weighted-bd-on-u}
			\sup_{x_n>0}x_n^{\frac{1}{q}+\frac{n-1}{p}}\|\overline{u}(\cdot,x_n)\|_{L^{\infty}(\rn)}\leq C\|f\|_{\dot{F}^{-\frac{1}{q}}_{p,q}(\rn)}.
		\end{equation}
		By the mean value property for the velocity field \cite[Theorem 4.5]{Sim} and with the same notation as above, we have
		\begin{align*}
			|\overline{u}_k(x)|&\leq \fint_{B_{x_n/2}(x)}|\overline{u}_k(y)|dy+\frac{1}{2}\fint_{B_{x_n/2}(x)}|\overline{\pi}(z)||z_k-x_k|dz:= I+II,\hspace{0.1cm}k=1,2,...,n.
		\end{align*}  
		Using Lemma \ref{lem:Aver}, H\"{o}lder's inequality and Proposition \ref{prop:Stokes-Triebel} simultaneously we estimate $I$ as follows. If $k=1,...,n-1$, then
		\begin{align}\label{u-1st-est}
			\nonumber I^q&\leq C|B_{\frac{x_n}{2}}(x)|^{-1} \bigg(\sum_{j=1}^{n-1}\int_{B_{x_n/2}(x)}|y_n^{-(n-1)}Q_{jk}(y_n^{-1}\cdot)\ast f_j|^qdy+\int_{B_{x_n/2}(x)}|(P_{k})_{y_n}f_n|^qdy\bigg)\\ 
			\nonumber&\leq Cx_n^{-n}\bigg(\sum_{j=1}^{n-1}\int_{E(B_{x_n}(x')\times [x_n/3,2x_n])}\iint_{\varGamma(z')}|y_n^{-(n-1-1/q)}Q_{jk}(y_n^{-1}\cdot)\ast f_j|^q\frac{dy'dy_n}{y_n^{n}}dz'+\\
			\nonumber   &\hspace{4.2cm}\int_{E(B_{x_n}(x')\times [x_n/3,2x_n])}\iint_{\varGamma(z')}|y_n^{1/q}(P_{k})_{y_n}f_n|^q\frac{dy'dy_n}{y_n^{n}}dz'\bigg)\\
			\nonumber&\leq Cx_n^{-1-\frac{(n-1)q}{p}}\bigg(\sum_{j=1}^{n-1}\big\|y_n^{-(n-1)}Q_{jk}(y_n^{-1}\cdot)\ast f_j\big\|^q_{T^{p,q}_{s_q}}+\big\|(P_{k})_{y_n}f_n\big\|^q_{T^{p,q}_{s_q}}\bigg)\\
			&\leq Cx_n^{-(1/q+(n-1)/p)q}\|f\|^q_{\dot{F}^{-1/q}_{p,q}(\rn)}.
		\end{align}
		When $k=n$, one may repeat the above steps with the kernels $P_j$ and $P_0$ respectively. 
		In order to estimate the integral $II:=\dfrac{1}{2}\displaystyle\fint_{B_{x_n/2}(x)}|\overline{\pi}(z)||z_k-x_k|dz$, we use the pressure estimate 
		\[\bigg(\fint_{B_{x_n}(x)}|\overline{\pi}(y)|^qdy\bigg)^{1/q}\leq cx_n^{-(1+1/q+(n-1)/p)}\|f\|_{\dot{F}^{-1/q}_{p,q}(\rn)},\quad \forall\, x\in \rnp\] 
		which can be derived from the steps leading to  \eqref{eq:ptw-pressure} and
		the  fact that if $z\in B_{x_n/2}(x)$, then $B_{x_n/2}(z)\subset B_{x_n}(x)$. Indeed, we have 
		\begin{align}\label{u-scd-est}
			\nonumber II&\leq |B_{x_n/2}(x)|^{-1}\int_{B_{x_n/2}(x)}\bigg(\fint_{B_{x_n/2}(z)}|\overline{\pi}(y)|dy\bigg)|z_k-x_k|dz\\
			\nonumber 	&\leq C|B_{x_n}(x)|^{-1}\bigg(\fint_{B_{x_n}(x)}|\overline{\pi}(y)|^qdy\bigg)^{1/q}\int_{B_{x_n/2}(x)}|z-x|dz\\
			\nonumber &\leq C|B_{x_n}(x)|^{-1}\bigg(\fint_{B_{x_n}(x)}|\overline{\pi}(y)|^qdy\bigg)^{1/q}\int_{0}^{x_n/2}r^{n}dr\\
			&\leq x_n^{-1/q-(n-1)/p}\|f\|_{\dot{F}^{-1/q}_{p,q}(\rn)}
		\end{align} 
		Combining \eqref{u-1st-est} and \eqref{u-scd-est}, we obtain \eqref{weighted-bd-on-u}. This achieves the proof of Lemma \ref{lem:lin-est}.
	\end{proof}  
	\subsection{Inhomogeneous Stokes system}
	Consider the operators $\mathscr{G}$ and $\Psi$ in $\rnp$ respectively defined by 
	\begin{align*}\mathscr{G}(F,H)(x)&=\int_{\rnp}G(x,y) F(y)dy+\int_{\rnp}\nabla_{y}  G(x,y) H(y)dy,\\
		\Psi(F,H)(x)&=\int_{\rnp}\g(x,y)F(y)dy+\int_{\rnp}\nabla_{y} \g(x,y) H(y)dy
	\end{align*}
	whenever the integrals make sense for almost every $x\in \rnp$. The kernels $G(x,y)=(G_{ij}(x,y))^{n}_{i,j=1}$ and $\g(x,y)=(g_j(x,y))^{n}_{j=1}$, $x\neq y$ are the Green tensors for the Stokes operator in $\rnp$, that is, coordinates-wise the functions satisfying
	\begin{align}\label{eq:Green-funct}
		\begin{cases}
			-\Delta_{x} G_{ij} +\partial_i g_j=\delta_{x}\delta_{ij}\hspace{0.2cm}\mbox{in}\hspace{0.2cm}\rnp\\
			\partial_{i} G_{ij}=0\hspace{0.2cm}\mbox{in}\hspace{0.2cm}\rnp\\
			G_{ij}(x,\cdot){\big|}_{\{(x',x_n):\hspace{0.1cm}x_n=0\}}=0
		\end{cases}
	\end{align} 
	in the sense of distributions where $\delta_{x}$ is the Dirac distribution with mass at $x\in \rnp$.
	Under mild assumptions on $F$ and $H$, the vector-valued functions $v=\mathscr{G}(F,H)$ and $w=\Psi (F,H)$ satisfy the system of equations 
	\begin{align}\label{eq:Inhom-Stokes}
		\begin{cases}
			-\Delta v+\nabla w=F-\dv H\hspace{0.2cm}\mbox{in}\hspace{0.2cm}\rnp\\
			\dv v=0\hspace{0.2cm}\mbox{in}\hspace{0.2cm}\rnp\\
			v=0\hspace{0.2cm}\mbox{on}\hspace{0.2cm}\partial\rnp.
		\end{cases}
	\end{align}
	Refined properties of Green matrices were recently obtained by the authors in \cite{KMT} relying on ideas introduced earlier in the articles \cite{Mazya} (for $n=2,3$) and \cite{Galdi} (for the general case). For our purpose we will need the following properties which include sharp pointwise decay bounds.  
	\begin{lemma}\label{lem:Green-pointwise-est}
		Let $n\geq 2$.	The Green tensor $G$ is symmetric, $G_{ij}(x,y)=G_{ji}(y,x)$ for all $x,y\in \rnp$, $x\neq y$ and satisfies together with  $\g$ the pointwise estimates
		\begin{align}
			\label{pt-bd}
			|G_{ij}(x,y)|\leq C\bigg(\dfrac{x_ny_n}{|x-y|^n}+\mathbf{1}_{\{n=2\}}\log(2+y_n|x-y|^{-1})\bigg)		
		\end{align}
		\begin{align}
			\label{grad-pt-bd}
			\big|\nabla^{\alpha}_{x}\nabla^{\beta}_{y}G_{ij}(x,y)\big|\leq C_{N}\begin{cases}|x-y|^{-(n-2+N)}\\
				\dfrac{x_ny_n}{|x-y|^{n+N}}\hspace{0.62cm}\mbox{if}\hspace{0.2cm}\alpha_n= \beta_{n}=0\\
				\dfrac{x_n}{|x-y|^{n-1+N}}\hspace{0.2cm}\mbox{if}\hspace{0.2cm}\alpha_n= 0
			\end{cases}	
		\end{align}
		for all multi-indices $\alpha,\beta$ with $|\alpha|+|\beta|=N>0$. Moreover,
		\begin{align}
			\label{g-bd}\big|\nabla^{\alpha}\g_j(x,y)\big|\leq C_{\alpha}|x-y|^{-(n-1)-|\alpha|},\hspace{0.1cm}j=1,...,n
		\end{align}
		where the constants are all independent of $x$ and $y$.
	\end{lemma}  
	These inequalities find their applicability in our next result which deals with the mapping properties of the potentials $\mathscr{G}$ and $\Psi$. Recall the spaces $\mathbf{X}_{p,q}$ and $\mathbf{Z}_{p,q}$ introduced in  Section \ref{s:1}.
	\begin{proposition}\label{prop:Green-pot-bound} Fix $n\geq 3$, let $p,q\in (1,\infty)$ with $q<p$; let $\sigma\in [1,q)$ and $\eta\in (1,q)$. Take $\tau\in (\eta,p)$ and $\va\in (\sigma,p)$ with the condition \[\frac{1}{\eta}+\frac{n-1}{\tau}=1+\frac{1}{\sigma}+\frac{n-1}{\va}=2+\frac{1}{q}+\frac{n-1}{p}.\] Then for all $F\in \X_{\tau,\eta}$ and $H\in \X_{\va,\sigma}$ we have $\mathscr{G}(F,H)\in \X_{p,q}$, $\Psi (F,H)\in \Z_{p,q}$ and  it holds that
		\begin{equation}\label{eq:nonlin-bound}
			\|\mathscr{G}(F,H)\|_{\X_{p,q}}+\|\Psi (F,H)\|_{\Z_{p,q}}\leq C(\|F\|_{\X_{\tau,\eta}}+\|H\|_{\X_{\varLambda,\sigma}})
		\end{equation}
		for some constant $C:=C(n,p,q)>0$ independent of $F$ and $H$.
	\end{proposition}
	\begin{remark}
		The proof of the above result reveals that elliptic estimates of the form 
		\begin{align}
			\sup_{x_n>0}x_n^{\frac{1}{q}+\frac{n-1}{p}+|\alpha|}\big\|\partial^{\alpha}u(\cdot,x_n)\big\|_{L^{\infty}(\rn)}\leq C(\|F\|_{\X_{\tau,\eta}}+\|H\|_{\X_{\varLambda,\sigma}})	
		\end{align} 
		are valid for each multi-index $\alpha$ with $\alpha_n=0$ where $u$ is a solution of the Stokes system \eqref{eq:Inhom-Stokes}. However, it is not clear whether vertical derivatives of $u$ enjoy this property. In fact, we are relying heavily on the second and third bound in \eqref{grad-pt-bd} which seem to fail in the case $\alpha_n\neq 0$ or $\beta_n\neq0$, see \cite[Remark 2.6]{KMT}.  	
	\end{remark}
	The proof of the proposition partially relies on the mapping properties in mixed Lebesgue spaces of the operator $G_{\beta}$ defined for $0<\beta<n$ by
	\begin{align}\label{Gbeta}
		G_{\beta}F(y)=\int_{\rnp}\dfrac{F(z)dz}{|y-z|^{n-\beta}}    
	\end{align}
	whenever the integral exists for almost all $y\in \rnp$.  For $p,q\in [1, \infty]$, let us denote by $L^{p}L^q(\rnp)$ the mixed Lebesgue space of functions $F:\rnp\rightarrow \mathbb{R}$ with the property that $x'\mapsto F(x',\cdot)\in L^p(\rn)$ and $x_n\mapsto F(\cdot,x_n)\in L^q(\mathbb{R}_+)$ and equip it with the norm \[\|F\|_{L^pL^q(\rnp)}=\big\|\|F(\cdot,x_n)\|_{L^{q}(\mathbb{R}_+,dx_n)}\big\|_{L^p(\rn)}.\]
	\begin{lemma}\label{lem:G-beta}
		Let $\beta\in (0,n)$, $\eta\in [1,\infty)$ and assume that $p,q,\tau\in (1,\infty)$ are such that $\tau< p<\infty$ and
		\begin{equation}
			\frac{1}{q}\leq \frac{1}{\eta}<\beta+\frac{1}{q},\quad  \frac{n-1}{p}=\frac{n-1}{\tau}+\frac{1}{\eta}-\frac{1}{q}-\beta.  
		\end{equation}
		Then the operator $G_{\beta}$ is bounded from $L^{\tau}L^{\eta}(\rnp)$ into $L^{p}L^q(\rnp)$.    
	\end{lemma}
	Recall the Riesz potential $I_{\alpha}$ of order $\alpha\in (0,n-1)$, that is,  the convolution operator with the kernel $|x'|^{\alpha-(n-1)}$, $x'\in \rn \setminus \{0\}$. 
	
	\begin{proof} 
		Along the lines of the proof of \cite[Lemma 2.2]{DY}, take $F\in L^{\tau}L^{\eta}(\rnp)$ and let $\widetilde{F}$ be the zero extension of $F$ to $\rnn$. For $1\leq\eta<\infty$ and $1<\tau<\infty$  we have
		\begin{align*}
			\|G_{\beta}F\|_{L^pL^q(\rnp)}=\bigg\|\big\|G_{\beta}F(y',\cdot)\big\|_{L^q(\mathbb{R}_+)}\bigg\|_{L^p(\rn)}.    
		\end{align*}
		Let $x'\in \rn$ and	set $S(x',s)=(|x'|^2+s^2)^{-\frac{n-\beta}{2}}$, $s\in \mathbb{R}$. For $1\leq \theta<\infty$ such that $\frac{1}{q}+1=\frac{1}{\theta}+\frac{1}{\eta}$ we use Minkowski's inequality to arrive at
		\begin{align*}
			\big\|G_{\beta}F(y',\cdot)\big\|_{L^q(\mathbb{R}_+)}&\leq\bigg\|\int_{\rnp}\dfrac{|F(z',z_n)|dz'dz_n}{(|y'-z'|^2+|\cdot-z_n|^2)^{\frac{n-\beta}{2}}}\bigg\|_{L^q(\mathbb{R}_+)}\\
			&\leq\bigg\|\int_{\rn}\big[S(y'-z',|\cdot|)\ast |\widetilde{F}|(z',\cdot)\big](y_n)dy'\bigg\|_{L^q(\mathbb{R}_+,dy_n)} \\
			&\leq C\int_{\rn}\big\|[S(y'-z',|\cdot|)\ast |\widetilde{F}|](z',\cdot)\|_{L^q(\mathbb{R})}dz'\\
			&\leq C\int_{\rn}\big\|S(y'-z',\cdot)\big\|_{L^{\theta}(\mathbb{R}_{+})}\|F(z',\cdot)\|_{L^{\eta}(\mathbb{R}_+)}dz'\\
			&\leq C[I_{\beta+\frac{1}{\theta}-1}\|F(\cdot,y_n)\|_{L^{\eta}(\mathbb{R}_+,dy_n)}](y'),\quad y'\in \rn
		\end{align*}
		Thus, if $\frac{n-1}{p}=\frac{n-1}{\tau}-(\beta+\frac{1}{\theta}-1)$, then by the boundedness of $I_{\alpha}$ in Lebesgue spaces, we find that
		\begin{equation*}
			\|G_{\beta}F\|_{L^pL^q(\rnp)}\leq C\big\|I_{\beta-\frac{1}{\theta}-1}\|F(y',\cdot)\|_{L^{\eta}(\mathbb{R}_+)}\big\|_{L^p(\rn,dy')}\leq C\big\|F\big\|_{L^\tau L^{\eta}(\rnp)}.
		\end{equation*}
	\end{proof}
	\begin{remark}\label{rmk:Gb-weighted}
		In the sequel, we will need an  analogue of Lemma \ref{lem:G-beta} in weighted mixed Lebesgue spaces of the form  
		\begin{align}\label{weighted-Gbeta}
			\bigg\|\big\|G_{\beta}F(\cdot,y_n)\big\|_{L^{q}(\mathbb{R}_+,y_n^qdy_n)}\bigg\|_{L^p(\rn)}\leq C\bigg\|\big\|F(\cdot,y_n)\big\|_{L^{\eta}(\mathbb{R}_+,y_n^{b\eta}dy_n)}\bigg\|_{L^r(\rn)}    
		\end{align}
		for all functions $F$ such that $(x',x_n)\mapsto x_n^{b}F\in L^{r}L^{\eta}(\rnp)$. This estimate is true under the following set of conditions
		\begin{align}\label{cds}
			\begin{cases}
				2+\dfrac{1}{q}=(n-1)\bigg(\dfrac{1}{r}-\dfrac{1}{p}\bigg)+\dfrac{1}{\eta}+b-(\beta-1)\\
				1<r<p<\infty,\hspace{0.2cm} b\geq 1\\
				n>\beta+2+\frac{1}{q}-\frac{1}{\eta}-b.
			\end{cases}    
		\end{align}
		In fact, one may use the same strategy as before to prove \eqref{weighted-Gbeta}. If $a\geq 1$ and $\delta>1$ are such that $$\dfrac{1}{\delta}+a=(n-1)\bigg(\dfrac{1}{r}-\dfrac{1}{p}\bigg)-(\beta-1),$$ then using the weighted convolution inequality \cite[Theorem 1.2]{GFWZ} for $n=1$, we obtain   
		\begin{align*}
			\big\|G_{\beta}F(y',\cdot)\big\|_{L^q(\mathbb{R}_+,y_n^qdy_n)}&\leq C\int_{\rn}\|S(y'-z',|\cdot|)\|_{L^{\delta}(\mathbb{R}_{+},y_n^{a\delta}dy_n)}\|F(z',\cdot)\|_{L^{\eta}(\mathbb{R}_+,y_n^{b\eta})}dz'\\
			&\leq I_{\frac{1}{\delta}+a+\beta-1}\|F(\cdot,y_n)\|_{L^{\eta}(\mathbb{R}_+,y_n^{b\eta}dy_n)}(y'),\quad y'\in \rn.
		\end{align*}
		This, in conjunction with \eqref{cds} gives the desired bound after taking the $L^p$-norm on both sides of the inequality. 
	\end{remark}
	We are now ready to prove Proposition \ref{prop:Green-pot-bound} and we divide the proof in two steps.\\
	\begin{step1}\label{step1} \normalfont The bound 
		\begin{equation}\label{eq:est-part-1}
			\|\mathscr{G}(F,H)\|_{\X_{p,q}}\leq C(\|F\|_{\X_{\tau,\eta}}+\|H\|_{\X_{\va,\sigma}}).
		\end{equation}
		Let $1<\eta<\infty$, $1\leq \sigma<\infty$ and $1<\tau,\va<\infty$ such that $\frac{1}{\eta}+\frac{n-1}{\tau}=2+\frac{1}{q}+\frac{n-1}{p}=1+\frac{1}{\sigma}+\frac{n-1}{\va}$. Pick $F$ in $\X_{\tau,\eta}$ and $H\in \X_{\va,\sigma}$. We first prove that
		\begin{equation}\label{eq:nonl-est-1}
			\displaystyle\sup_{x_n>0}x_n^{1/q+(n-1)/p}\|\mathscr{G}(F,H)(\cdot,x_n)\|_{L^{\infty}(\rn)}\leq \|F\|_{\X_{\tau,\eta}}.
		\end{equation}
		Fix $x'\in \rn$ and $x_n>0$ and write
		\begin{equation*}
			\int_{\rnp}G(x',x_n,y) F(y)dy=J_1+J_2+J_3+J_4
		\end{equation*}
		where \[J_1=\int_{B_{x_n}(x')}\int_{0}^{x_n/2}G(x,y) F(y)dy,\hspace{0.1cm}J_2=\int_{B_{x_n}(x')}\int_{x_n/2}^{2x_n}G(x,y) F(y)dy,\]
		\[J_3=\int_{\rn\setminus B_{x_n}(x')}\int_{0}^{2x_n}G(x,y) F(y)dy,\hspace{0.1cm}J_4=\int_{\rn}\int_{2x_n}^{\infty}G(x,y) F(y)dy.\]
		Next, we estimate each of these integrals by means of the pointwise inequalities from Lemma \ref{lem:Green-pointwise-est}. Indeed, starting with $J_1$ and using the summation convention, we have   
		\begin{align*}
			|J_1|&\leq\int_{B_{x_n}(x')}\int_{0}^{x_n/2}|G(x',x_n,y)||F(y)|dy\\
			&\leq C\int_{B_{x_n}(x')}\int_{0}^{x_n/2}\frac{|F(y)|}{(|x'-y'|^{2}+(x_n-y_{n})^{2})^{\frac{n-2}{2}}}dy_{n}dy'\\
			&\leq Cx_n^{-(n-2)+\frac{n}{\eta'}}\bigg(\int_{B_{x_n}(x')}\int_{0}^{x_n/2}|F(y)|^{\eta}dy_{n}dy'\bigg)^{1/\eta}\\
			&\leq Cx_n^{-(n-2)+\frac{n}{\eta'}}\bigg(\int_{\rn}\iint_{\varGamma(z') \cap B_{x_n}(x')\times (0,x_n/2)}|F(y)|^{\eta}y_{n}^{-(n-1)}dy_{n}dy'dz'\bigg)^{1/\eta}\\	
			&\leq Cx_n^{-(n-2)+\frac{n}{\eta'}}\|F\|_{T^{\tau,\eta}_{s_{\eta}}}\big|E(B_{x_n/2}(x')\times (0,x_n/2))\big|^{\frac{\tau-\eta}{\tau\eta}}\\
			&\leq Cx_n^{2-\frac{n-1}{\tau}-\frac{1}{\eta}}\|F\|_{\X_{\tau,\eta}},
		\end{align*}
		where we have utilized H\"{o}lder's inequality in order to derive the third and fifth bounds in the above chain of estimates and $\frac{1}{\eta'}+\frac{1}{\eta}=1$. On the other hand, 
		\begin{align*}
			|J_2|&\leq \int_{B_{x_n}(x')}\int_{x_n/2}^{2x_n}|G(x,y)||F(y)|dy\\
			&\leq C\int_{B_{x_n}(x')}\int_{x_n/2}^{2x_n}|x-y|^{-(n-2)}|F(y)|dy\\
			&\leq C\sup_{y_{n}>0}y_n^{\frac{n-1}{\tau}+\frac{1}{\eta}}\|F(\cdot,y_{n})\|_{L^{\infty}(\rn)}\int_{B_{x_n}(x')}\int_{x_n/2}^{2x_n}\frac{y_{n}^{-\frac{1}{\eta}-\frac{n-1}{\tau}}dy_{n}dy'}{[|x'-y'|^{2}+(x_n-y_{n})^{2}]^{(n-2)/2}}\\
			&\leq Cx_n^{-\frac{1}{\eta}-\frac{n-1}{\tau}}\|F\|_{\X_{\tau,\eta}}\int_{B_{x_n}(x')}\int_{x_n/2}^{2x_n}|x'-y'|^{-(n-2)}dy_{n}dy'\\
			&\leq Cx_n^{1-\frac{1}{\eta}-\frac{n-1}{\tau}}\|F\|_{\X_{\tau,\eta}}\int_{B_{x_n}(x')}|x'-y'|^{-(n-2)}dy'\\
			&\leq Cx_n^{2-\frac{1}{\eta}-\frac{n-1}{\tau}}\|F\|_{\X_{\tau, \eta}}.
		\end{align*}
		Similarly as above, by Lemma \ref{lem:Aver} and H\"{o}lder's inequality, we find that
		\begin{align*}
			|J_3|&\leq \int_{\rn\setminus B_{x_n}(x')}\int_{0}^{2x_n}|G(x,y)||F(y)|dy\\
			&\leq Cx_n\int_{\rn\setminus B_{x_n}(x')}\int_{0}^{2x_n}|x-y|^{-(n-1)}|F(y)|dy\\
			&\leq Cx_n\sum_{k=1}^{\infty}\int_{2^{k}B_{x_n}(x')\setminus 2^{k-1}B_{x_n}(x')}\int_{0}^{2x_n}|x-y|^{-(n-1)}|F(y)|dy\\
			&\leq Cx_n^{2-n+n/\eta'}\sum_{k=1}^{\infty}2^{-(k-1)(n-1)+\frac{(n-1)k}{\eta'}}\bigg(\int_{2^{k}B_{x_n}(x')}\int_{0}^{2x_n}|F(y)|^{\eta}dy_{n}dy'\bigg)^{\frac{1}{\eta}}\\	
			&\leq Cx_n^{2-\frac{n}{\eta}+\frac{(n-1)(\tau-\eta)}{\tau\eta}}\bigg[\int_{\rn}\bigg(\iint_{\varGamma(z')}|F(y)|^{\eta}\frac{dy}{y_n^{n-1}}\bigg)^{\tau/\eta}dz'\bigg]^{\frac{1}{\tau}}\bigg(\sum_{k=1}^{\infty}2^{-\frac{k(n-1)}{\tau}}\bigg)\\
			&\leq Cx_n^{2-\frac{1}{\eta}-\frac{n-1}{\tau}}\|F\|_{\X_{\tau,\eta}}.
		\end{align*}
		Again, by using the Green matrix bound \eqref{pt-bd},  we bound $J_4$ as follows 
		\begin{align*}
			|J_4|&\leq \int_{\rn}\int_{2x_n}^{\infty}|G(x,y)||F(y)|dy\\
			&\leq C\int_{\rn}\int_{2x_n}^{\infty}\frac{x_ny_n|F(y)|}{|x-y|^{n}}dy\\
			&\leq C\int_{\rn}\int_{2x_n}^{\infty}\frac{x_ny_n|F(y)|dy_{n}dy'}{\big[|x'-y'|^{2}+y_{n}^{2}\big]^{\frac{n}{2}}}\\
			&\leq C\sup_{x_n>0}x^{\frac{1}{\eta}+\frac{n-1}{\tau}}_n\|F(\cdot,x_n)\|_{L^{\infty}(\rn)}\bigg(\int_{2x_n}^{\infty}x_ny^{-\frac{1}{\eta}-\frac{n-1}{\tau}}_ndy_n\bigg)\bigg(\int_{\rn}\frac{dz'}{\big[|z'|^{2}+1\big]^{\frac{n}{2}}}\bigg)\\
			&\leq Cx_n^{2-\frac{1}{\eta}-\frac{n-1}{\tau}}\|F\|_{\X_{\tau,\eta}}.
		\end{align*}
		In the same vein, we establish the weighted gradient sup-norm estimate 
		\begin{equation}\label{eq:weighted-sup-est}
			\sup_{x_n>0}x_n^{1/q+p/(n-1)}\bigg\|\int_{\rnp}\nabla_{y} G(\cdot ,x_n,y)H(y)dy\bigg\|_{L^{\infty}(\rn)}\leq C\|H\|_{\X_{\va,\sigma}}.
		\end{equation}
		Decompose the solid integral in the above estimate into four parts to get
		\[L_1=\int_{B_{x_n}(x')}\int_{0}^{x_n/2}\nabla_{y} G(x,y) H(y)dy,\hspace{0.1cm}L_2=\int_{B_{x_n}(x')}\int_{x_n/2}^{2x_n}\nabla_{y} G(x,y) H(y)dy,\]
		\[L_3=\int_{\rn\setminus B_{x_n}(x')}\int_{0}^{2x_n}\nabla_{y} G(x,y) H(y)dy,\hspace{0.1cm}L_4=\int_{\rn}\int_{2x_n}^{\infty}\nabla_{y} G(x,y) H(y)dy.\]
		Utilizing \eqref{grad-pt-bd},  H\"{o}lder's  inequality, Lemmas \ref{lem:Aver} and \ref{lem:E_S} we arrive at   
		\begin{align*}
			|L_1|&\leq\int_{B_{x_n}(x')}\int_{0}^{x_n/2}|\nabla_{y} G(x,y)||H(y)|dy\\
			&\leq C\int_{B_{x_n}(x')}\int_{0}^{x_n/2}\frac{|H(y',y_{n})|}{(|x'-y'|^{2}+(x_n-y_{n})^{2})^{\frac{n-1}{2}}}dy_{n}dy'\\
			&\leq Cx_n^{1-\frac{1}{\sigma}-\frac{n-1}{\va}}\|\mathcal{A}_{\sigma}(y_n^{1/\sigma}H)\|_{L^{\va}(\rn)}\\	
			&\leq Cx_n^{1-\frac{1}{\sigma}-\frac{n-1}{\va}}\|H\|_{\X_{\va,\sigma}}.
		\end{align*}
		Next, noticing that $|\nabla G_{ij}(x,\cdot)|$ belongs to the weak-Lebesgue space $L^{\frac{n}{n-1},\infty}(\rnp)$ uniformly for all $x\in \rnp$, it follows that
		\begin{align*}
			|L_2|&\leq \int_{B_{x_n}(x')}\int_{x_n/2}^{2x_n}|\nabla_{y} G(x,y) H(y)|dy\\
			&\leq C\sup_{x_{n}>0}x_{n}^{\frac{1}{\sigma}+\frac{n-1}{\va}}\|H(\cdot,x_{n})\|_{L^{\infty}(\rn)}\int_{B_{x_n}(x')}\int_{x_n/2}^{2x_n}y_{n}^{-\frac{1}{\sigma}-\frac{n-1}{\va}}|\nabla_{y} G(x,y)|dy_{n}dy'\\
			&\leq Cx_n^{-\frac{1}{\sigma}-\frac{n-1}{\va}}\sup_{x_{n}>0}x_{n}^{\frac{1}{\sigma}+\frac{n-1}{\va}}\|H(\cdot,x_{n})\|_{L^{\infty}(\rn)}\|\nabla_{y} G(x,\cdot)\|_{L^{1}(B_{x_n}(x')\times [x_n/2,2x_n])}\\
			&\leq Cx_n^{1-\frac{1}{\sigma}-\frac{n-1}{\va}}\|H\|_{\X_{\va,\sigma}}.
		\end{align*}
		Recall here that for any $p>1$ the belonging of $f$ to $L^{p,\infty}(\rn)$ is equivalent to the condition
		\begin{align*}
			\sup_{E\subset \rn}|E|^{1/p-1}\int_{E}|f(y)|dy<\infty    
		\end{align*}
		where the supremum runs over all open set $E$ of $\rn$. 
		We argue as before to bound $L_3$ 
		\begin{align*}
			|L_3|&\leq \int_{\rn\setminus B_{x_n}(x')}\int_{0}^{2x_n}|\nabla_{y} G(x,y)H(y)|dy_{n}dy'\\
			&\leq \sum_{k=1}^{\infty}\int_{2^{k}B_{x_n}(x')\setminus 2^{k-1}B_{x_n}(x')}\int_{0}^{2x_n}|\nabla_{y} G(x,y)||H(y)|dy_{n}dy'\\
			&\leq C\sum_{k=1}^{\infty}\int_{2^{k}B_{x_n}(x')\setminus 2^{k-1}B_{x_n}(x')}\int_{0}^{2x_n}|x-y|^{-n+1}|H(y)|dy_{n}dy'\\
			&\leq Cx_n^{1-\frac{n}{\sigma}}\sum_{k=1}^{\infty}2^{-(n-1)k+\frac{k(n-1)}{\sigma'}}\bigg(\int_{\rn}\iint_{\varGamma(z')\cap [2^kB_{x_n}(x')\times (0,2x_n)]}|H(y)|^{\sigma}\frac{dy_{n}dy'}{y_{n}^{n-1}}dz'\bigg)^{\frac{1}{\sigma}}\\
			&\leq Cx_n^{-(n-1)+\frac{n}{\sigma'}+\frac{(\varLambda-\sigma)}{\varLambda\sigma}n}\|\mathcal{A}_{\sigma}(y_n^{1/\sigma}H)\|_{L^{\varLambda}(\rn)}\bigg(\sum_{k=1}^{\infty}2^{-(k-1)(n-1)+\frac{k(n-1)}{\sigma'}+k(n-1)\frac{(\varLambda-\sigma)}{\sigma\varLambda}}\bigg)\\
			&\leq Cx_n^{1-\frac{1}{\sigma}-\frac{n-1}{\va}}\|H\|_{T^{\va,\sigma}_{s_{\sigma}}}\sum_{k=1}^{\infty}2^{-\frac{(n-1)}{\va}k}\\
			&\leq Cx_n^{1-\frac{1}{\sigma}-\frac{n-1}{\va}}\|H\|_{\X_{\va,\sigma}}.
		\end{align*}
		Finally, observe that for $y_{n}>2x_n$, we have $y_{n}-x_n>\frac{1}{2}y_{n}$ so that by the third bound in  \eqref{grad-pt-bd}, we find that
		\begin{align*}
			|L_4|&\leq \int_{\rn}\int_{2x_n}^{\infty}|\nabla_{y} G(x,y)||H(y)|dy\\
			&\leq C\int_{\rn}\int_{2x_n}^{\infty}\frac{x_n|H(y)|dy_{n}dy'}{\big[|x'-y'|^{2}+(x_n-y_{n})^{2}\big]^{\frac{n}{2}}}\\
			&\leq C\int_{\rn}\int_{2x_n}^{\infty}\frac{x_n|H(y)|dy_{n}dy'}{\big[|x'-y'|^{2}+y_{n}^{2}\big]^{\frac{n}{2}}}\\
			&\leq C\sup_{x_n>0}x_n^{\frac{1}{\sigma}+\frac{n-1}{\varLambda}} \|H(\cdot,x_n)\|_{L^{\infty}(\rn)}\bigg(\int_{\rn}\frac{dy'}{(|y'|^{2}+1)^{\frac{n}{2}}}\bigg)\bigg(\int_{2x_n}^{\infty}x_ny^{-\frac{1}{\sigma}-\frac{n-1}{\varLambda}-1}_ndy_n\bigg)\\
			&\leq Cx_n^{1-\frac{1}{\sigma}-\frac{n-1}{\va}}\|H\|_{\X_{\va,\sigma}}.
		\end{align*}
		Summing up  all the above inequalities, one obtains \eqref{eq:weighted-sup-est}. Next, we show that 
		\begin{align}\label{tent-est}
			\|\mathscr{G}(F,H)\|_{T^{p,q}_{s_{q}}}\leq C(\|F\|_{T^{\tau,\eta}_{s_{\eta}}}+\|H\|_{T^{\va,\sigma}_{s_{\sigma}}}).
		\end{align}  
		Write
		\begin{align*}
			\|\mathscr{G}(F,H)\|_{T^{p,q}_{s_{q}}}&\leq \bigg\|\int_{\rnp}G(\cdot,y) F(y)dy\bigg\|_{T^{p,q}_{s_{q}}}+\bigg\|\int_{\rnp}\nabla_{y}  G(\cdot,y) H(y)dy\bigg\|_{T^{p,q}_{s_{q}}}\\
			&:=I+II.    
		\end{align*}
		Fix $x'\in \rn$ and $y_n>0$ and let's decompose $F\in L^{\eta}_{loc}(\rnp)$ into three parts \[F=F\textbf{1}_{B_{4y_n}(x')\times (0,4y_n]}+F\textbf{1}_{B_{4y_n}(x')\times (4y_n,\infty)}+F\textbf{1}_{(\rn\setminus B_{4y_n}(x'))\times (0,\infty)}=F^1+F^2+F^3\]
		and write
		\begin{align*}
			I:=\Sigma_1+\Sigma_2+\Sigma_3,\quad \Sigma_i=\bigg\|\int_{\rnp}G(\cdot,y) F^i(y)dy\bigg\|_{T^{p,q}_{s_{q}}},\quad i=1,2,3.    
		\end{align*}
		We control  $\Sigma_3$ using the following
		\begin{claim}\label{claim:1}
			For all $x'\in \rn$ and $y_n>0$, there exists $C>0$ independent on $x'$, $y_n$ and $F$ such that
			\begin{align*}
				A(x',y_n)&\leq CG_{2}\bigg(\fint_{B_{y_n}(\cdot)}|F(z',\cdot)|dz'\bigg)(x',y_n).
			\end{align*}
			Here,  
			\begin{align*}
				A(x',y_n)=\bigg(\fint_{B_{y_n}(x')}\bigg|\int_{\rnp}G(y,z) F^3(z)dz\bigg|^qdy'\bigg)^{\frac{1}{q}}, \quad (x',y_n)\in \rnp.    
			\end{align*}
		\end{claim}
		\begin{proof}
			We have
			\begin{align*}
				A(x',y_n)&\leq \bigg(\fint_{B_{y_n}(x')}\bigg(\int_{\rnp}|G(y,z)||F^3(z)|dz\bigg)^{q}dy'\bigg)^{1/q}\\
				&\leq C\bigg(\fint_{B_{y_n}(x')}\bigg(\int_{0}^{\infty}\int_{\rn\setminus B_{4y_n}(x')}\dfrac{|F(z',z_n)|dz'dz_n}{(|y'-z'|^2+|y_n-z_n|^2)^{\frac{n-2}{2}}}\bigg)^{q}dy'\bigg)^{1/q}\\
				&\leq C\bigg(\fint_{B_{y_n}(x')}\bigg(\int_{0}^{\infty}\int_{\{|x'-z'|>4y_n\}}\dfrac{|F(z',z_n)|dz'dz_n}{(|y'-z'|^2+|y_n-z_n|^2)^{\frac{n-2}{2}}}\fint_{B_{y_n}(z')}dw\bigg)^{q}dy'\bigg)^{\frac{1}{q}}\\
				&\leq C\bigg(\fint_{B_{y_n}(x')}\bigg(\int_{0}^{\infty}\int_{\{|x'-w|>3y_n\}}\fint_{B_{y_n}(w)}\dfrac{|F(z',z_n)|dz'dwdz_n}{(|y'-z'|^2+|y_n-z_n|^2)^{\frac{n-2}{2}}}\bigg)^{q}dy'\bigg)^{1/q}\\
				&\leq C\bigg(\fint_{B_{y_n}(x')}\bigg(\int_{0}^{\infty}\int_{\{|x'-w|>3y_n\}}(|x'-w|^2+|y_n-z_n|^2)^{\frac{-(n-2)}{2}}\cdot\\
				&\hspace{2.4cm}\fint_{B_{y_n}(w)}\bigg[\dfrac{|x'-w|^2+|y_n-z_n|^2}{|y'-z'|^2+|y_n-z_n|^2}\bigg]^{\frac{n-2}{2}}|F(z',z_n)|dz'dwdz_n\bigg)^{q}dy'\bigg)^{\frac{1}{q}}\\
				&\leq C\bigg(\fint_{B_{y_n}(x')}\bigg(\int_{0}^{\infty}\int_{\{|x'-w|>3y_n\}}(|x'-w|^2+|y_n-z_n|^2)^{\frac{-(n-2)}{2}}\cdot\\
				&\hspace{7cm}\fint_{B_{y_n}(w)}|F(z',z_n)|dz'dwdz_n\bigg)^{q}dy'\bigg)^{1/q}\\
				&\leq C\int_{0}^{\infty}\int_{\{|x'-w|>3y_n\}}(|x'-w|^2+|y_n-z_n|^2)^{\frac{-(n-2)}{2}}\bigg(\fint_{B_{y_n}(w)}|F(z',z_n)|dz'\bigg)dwdz_n\\
				&\leq C\int_{\rnp}(|x'-w|^2+|y_n-z_n|^2)^{\frac{-(n-2)}{2}}\bigg(\fint_{B_{y_n}(w)}|F(z',z_n)|dz'\bigg)dwdz_n\\
				&\leq CG_{2}\bigg(\fint_{B_{y_n}(\cdot)}|F(z',\cdot)|dz'\bigg)(x',y_n).
			\end{align*}
		\end{proof}
		Applying Lemma \ref{lem:G-beta} and Jensen's inequality, the above claim clearly implies that 
		\begin{align*}
			\Sigma_3= \|A\|_{L^pL^q(\rnp)}&\leq C\bigg\|G_2\bigg(\fint_{B_{y_n}(\cdot)}|F(z',\cdot)|dz'\bigg)\bigg\|_{L^pL^q(\rnp)}\\
			&\leq C\bigg\|\fint_{B_{y_n}(\cdot)}|F(z',y_n)|dz'\bigg\|_{L^{\tau}L^{\eta}(\rnp)}\leq C\|F\|_{T^{\tau,\eta}_{s_{\eta}}}.
		\end{align*}
		To bound $\Sigma_2$, we first observe that 
		\begin{align}\label{I^2}
			\nonumber\bigg|\int_{\rnp}G(y,z) F^2(z)dz\bigg|&\leq Cy_n\big[\mathcal{A}^4_{\eta}(z_n^{1/\eta}F)\big](x')\bigg(\int_{4y_n}^{\infty}\int_{B_{4y_n}(x')}\dfrac{z_n^{1+\frac{n}{\eta-1}}dz'dz_n}{(|y'-z'|^2+|y_n-z_n|^2)^{\frac{n\eta'}{2}}}\bigg)^{\frac{1}{\eta'}}\\
			&\leq Cy_n^{2-\frac{1}{\eta}}\big[\mathcal{A}^4_{\eta}(z_n^{1/\eta}F)\big](x'),\quad x'\in B(y',y_n).
		\end{align}
		On the other hand, this inequality also implies, thanks to Remark \ref{rmk:ind-aperture} the pointwise bound
		\begin{align}\label{I^22}
			\bigg|\int_{\rnp}G(y,z) F^2(z)dz\bigg|&\leq Cy_{n}^{-\frac{n-1}{\tau}+2-\frac{1}{\eta}}\|\mathcal{A}_{\eta}(z_n^{1/\eta}F)\|_{L^{\tau}(\rn)},\quad  y'\in B(x',y_n).
		\end{align}
		Let $M>0$ to be determined later. Using \eqref{I^2} and \eqref{I^22}, we find that 
		\begin{align*}
			\iint_{\varGamma(x')}\bigg|y_n^{1/q}\int_{\rnp}G(y,z) F^2(z)dz\bigg|^q\frac{dy}{y_n^{n}}&\leq \int_{0}^{M}\fint_{B_{y_n}(x')}\bigg|\int_{\rnp}G(y,z) F^2(z)dz\bigg|^qdy'dy_n+\\
			&\hspace{1cm}\int_{M}^{\infty}\fint_{B_{y_n}(x')}\bigg|\int_{\rnp}G(y,z) F^2(z)dz\bigg|^qdy'dy_n\\
			&\leq CM^{1+(2-\frac{1}{\eta})q}[\mathcal{A}^4_{\eta}(z_n^{\frac{1}{\eta}}F)(x')]^q+M^{-\frac{(n-1)q}{p}}\|F\|^q_{T^{\tau,\eta}_{s_{\eta}}}. 
		\end{align*}
		Optimizing this inequality with respect to $M$, that is taking $M=\bigg(\dfrac{\|F\|_{T^{\tau,\eta}_{s_{\eta}}}}{\big[\mathcal{A}^4_{\eta}(z_n^{1/\eta}F)\big](x')}\bigg)^{\frac{\tau}{n-1}}$, we arrive at
		\begin{align*}
			\bigg(\iint_{\varGamma(x')}\bigg|y_n^{\frac{1}{q}}\int_{\rnp}G(y,z) F^2(z)dz\bigg|^qy_n^{-n}dy'dy_n\bigg)^{\frac{1}{q}}&\leq C\|F\|^{1-\frac{\tau}{p}}_{T^{\tau,\eta}_{s_{\eta}}}\big[\mathcal{A}^4_{\eta}(z_n^{1/\eta}F)(x')]^{\frac{\tau}{p}},\, x'\in \rn.
		\end{align*}
		Taking the $L^{p}$-norm on both sides of the inequality and using Remark \ref{rmk:ind-aperture}, we conclude that
		\begin{align*}
			\Sigma_2\leq C\|F\|_{T^{\tau,\eta}_{s_{\eta}}}.    
		\end{align*}
		Finally, estimating $\Sigma_1$ goes through a duality argument. Let $\rho=p/q$ and $\varphi\in L^{\rho'}(\rn)$, $\varphi\geq 0$ and define the operator $M_t$, $t>0$ by  \[M_{t}\varphi(x')=t^{-(n-1)}\int_{B_t(x')}\varphi(y')dy',\quad t>0.\] If $\langle\cdot,\cdot\rangle$ denotes the duality bracket between $L^p(\rn)$ and its dual $L^{p'}(\rn)$, then
		\begin{align*}
			\bigg\langle \mathcal{A}^q_q\bigg[y_n^{\frac{1}{q}}\int_{\rnp}G(\cdot,z) F^1(z)dz\bigg],\varphi\bigg\rangle &= \int_{\rn}\int_{\varGamma(x')}\bigg|\int_{\rnp}G(y,z) F^1(z)dz\bigg|^q\dfrac{dy'dy_n}{y_n^{n-1}}\varphi(x')dx'\\
			&\leq C\int_{\rn}\int_0^{\infty}\fint_{B_{y_n}(y')}\varphi(x')dx'[G_2|F|(y',y_n)]^qdy'dy_n\\
			&\leq C\int_{\rn}\int_0^{\infty}[G_2|F|(y',y_n)]^qM_{y_n}\varphi(y')dy_ndy'\\
			&\leq C\|G_2F\|^q_{L^{p}L^q(\rnp)}\|M_{\cdot}\varphi\|_{L^{\rho'}L^{\infty}(\rnp)}\\
			&\leq C\|G_2F\|^q_{L^{p}L^q(\rnp)}\|\mathcal{M}\varphi\|_{L^{\rho'}(\rn)}.
		\end{align*}
		The following inequality is a consequence of Lebesgue's differentiation Theorem and Fatou's lemma:		\[\int_0^{\infty}|F(y',y_n)|^{\eta}dy_n\leq c_n\liminf_{\alpha\rightarrow 0}[\mathcal{A}^{\alpha}_{\eta}(y_n^{1/\eta}F)(y')]^{\eta},\,\,\,y'\in \rn.\]
		Applying Lemma \ref{lem:G-beta}, the boundedness of the Hardy-Littlewood maximal function $\mathcal{M}$ in Lebesgue spaces successively, we obtain 
		\begin{align*}
			\bigg\langle \mathcal{A}^q_q\bigg[y_n^{1/q}\int_{\rnp}G(\cdot,z) F^1(z)dz\bigg],\varphi\bigg\rangle &\leq C\|F\|^q_{T^{\tau,\eta}_{s_{\eta}}}\|\varphi\|_{L^{\rho'}(\rn)}\quad \forall \,\varphi \in L^{\rho'}(\rn),    
		\end{align*}
		from which it plainly follows that 
		\begin{align*}
			\Sigma_1\leq C    \|F\|_{T^{\tau,\eta}_{s_{\eta}}}.
		\end{align*}
		We equally estimate $II$ splitting $H$ into three components exactly as before and follow the same procedure (details are left to the interested reader). This yields 
		\begin{align*}
			\bigg\|\int_{\rnp}\nabla_yG(\cdot,y) H(y)dy\bigg\|_{T^{p,q}_{s_{q}}}\leq C\|H\|_{T^{\va,\sigma}_{s_{\sigma}}}.
		\end{align*}
		Summarizing, we see that \eqref{eq:est-part-1} holds true. This finishes Step 1.
	\end{step1} 
	\textbf{Step 2.} The estimate 
	\begin{align}\label{eq:est-part-2}
		\|\Psi (F,H)\|_{\Z_{p,q}}\leq C\big(\|F\|_{\X_{\tau,\eta}}+\|H\|_{\X_{\va,\sigma}}\big)
	\end{align}
	for all $F\in \X_{\tau,\eta}$ and $H\in \X_{\va,\sigma}$. 
	We have 
	\begin{align*}
		\|\Psi(F,H)\|_{T_{\sigma_q}^{p,q}}\leq \bigg\|\int_{\rnp}\g(\cdot,z) F(z)dz\bigg\|_{T_{\sigma_q}^{p,q}}+\bigg\|\int_{\rnp}\nabla_{z}  \g(\cdot,z) H(z)dz\bigg\|_{T^{p,q}_{\sigma_q}}:=III+IV.    
	\end{align*}
	Let $F^1$, $F^2$ and $F^3$ as above and write correspondingly
	\begin{align*}
		III\leq III_1+III_2+III_3,\quad III_i=\bigg\|\int_{\rnp}\g(\cdot,z) F^i(z)dz\bigg\|_{T_{\sigma_q}^{p,q}},\quad i=1,2,3.    
	\end{align*}
	From the proof of Claim \ref{claim:1}, we easily obtain the estimate
	\begin{align*}
		\bigg(\fint_{B_{y_n}(x')}\bigg|\int_{\rnp}\g(y,z) F^3(z)dz\bigg|^qdy'\bigg)^{\frac{1}{q}}\leq CG_{1}\bigg(\fint_{B_{y_n}(\cdot)}|F(z',\cdot)|dz'\bigg)(x',y_n),\quad (x',y_n)\in \rnp.
	\end{align*}
	Now let $r\in (\tau,p)$ such that $\frac{1}{r}+\frac{1}{n-1}\leq \frac{1}{\tau}$.  Invoking \eqref{weighted-Gbeta} with $\beta=1$ and  $b=(n-1)(\frac{1}{\tau}-\frac{1}{r})$ together with Jensen's inequality we arrive at
	\begin{align*}
		III_3&\leq C\bigg\|\big\|G_1\bigg(\fint_{B_{y_n}(\cdot)}|F(z',\cdot)|dz'\bigg)\big\|_{L^q(\mathbb{R}_+,y^q_ndy_n)}\bigg\|_{L^p(\rn)}\\
		&\leq C\bigg(\int_{\rn}\bigg(\int_{0}^{\infty}\fint_{B_{y_n}(x')}|y_n^{b}F(z',y_n)|^{\eta}dz'dy_n\bigg)^{\frac{r}{\eta}}dx'\bigg)^{\frac{1}{r}}\\
		&\leq C	\big\|\mathcal{A}_{\eta}(y_n^{b+1/\eta}F)\big\|_{L^r(\rn)}\leq C\|F\|_{T^{\tau,\eta}_{s_{\eta}}}.
	\end{align*}
	The last inequality follows from the embedding \eqref{embedding-tents} (with $s_1=-\frac{1}{\eta(n-1)}$, $s_2=-\frac{b+1/\eta}{n-1}$, $q=\eta$, $p_1=\tau$ and $p_2=r$). Moving on, we use \eqref{g-bd} and H\"{o}lder's inequality to get the pointwise bound 
	\begin{align}\label{III}
		\nonumber\bigg|\int_{\rnp}\g(y,z) F^2(z)dz\bigg|&\leq C|G_1F^2(y)|\\
		\nonumber&\leq C[\mathcal{A}^4_{\eta}(y_n^{1/\eta}F)](x')\bigg(\int_{4y_n}^{\infty}\int_{B_{4y_n}(x')}\dfrac{z_n^{\frac{n-1}{\eta-1}}dz'dz_n}{(|y'-z'|^2+|y_n-z_n|^2)^{\frac{(n-1)\eta'}{2}}}\bigg)^{\frac{1}{\eta'}}\\
		&\leq Cy_n^{1-\frac{1}{\eta}}\big[\mathcal{A}^4_{\eta}(z_n^{1/\eta}F)\big](x'),\quad x'\in B_{y_n}(y')
	\end{align}
	from which it follows that
	\begin{align}\label{III2}
		\bigg|\int_{\rnp}\g(y,z) F^2(z)dz\bigg|&\leq Cy_{n}^{-\frac{n-1}{\tau}-\frac{1}{\eta}+1}\|F\|_{T^{\tau,\eta}_{s_{\eta}}},\quad  y'\in B(x',y_n).
	\end{align}
	Therefore, for $M>0$ to be determined later, we have
	\begin{align*}
		\iint_{\varGamma(x')}y_n^{q+1}\bigg|\int_{\rnp}\g(y,z) F^2(z)dz\bigg|^q\frac{dy'dy_n}{y_n^n}
		&\leq \int_{0}^{M}\fint_{B_{y_n}(x')}y_n^q\bigg|\int_{\rnp}\g(y,z) F^2(z)dz\bigg|^qdy'dy_n+\\
		&\hspace{0.41cm}\int_{M}^{\infty}\fint_{B_{y_n}(x')}y_n^q\bigg|\int_{\rnp}\g(y,z)F^2(z)dz\bigg|^qdy'dy_n\\
		&\lesssim M^{1+(2-\frac{1}{\eta})q}\big[\mathcal{A}^4_{\eta}(z_n^{\frac{1}{\eta}}F)(x')\big]^q+M^{-\frac{n-1}{p}q}\|F\|^q_{T^{\tau,\eta}_{s_{\eta}}}. 
	\end{align*}
	The choice $M=\bigg(\dfrac{\|F\|_{T^{\tau,\eta}_{s_{\eta}}}}{\big[\mathcal{A}^4_{\eta}(z_n^{1/\eta}F)(x')\big]}\bigg)^{\frac{\tau}{n-1}}$ yields the bound
	\begin{align*}
		\bigg(\iint_{\varGamma(x')}y_n^{1+q}\bigg|\int_{\rnp}\g(y,z) F^2(z)dz\bigg|^q\frac{dy'dy_n}{y_n^n}\bigg)^{\frac{1}{q}}&\leq C\|F\|^{1-\frac{\tau}{p}}_{T^{\tau,\eta}_{s_{\eta}}}\big[\mathcal{A}^4_{\eta}(z_n^{1/\eta}F)(x')\big]^{\frac{\tau }{p}},\quad x'\in \rn.
	\end{align*}
	Hence, (after taking the $L^{p}$-norm on both sides of the previous inequality)
	\begin{align*}
		III_2\leq C\|F\|_{T^{\tau,\eta}_{s_{\eta}}}.    
	\end{align*}
	We also claim that $III_1\leq C    \|F\|_{T^{\tau,\eta}_{s_{\eta}}}$.
	In fact, setting $\displaystyle VF(y',y_n)=y_n^{1+\frac{1}{q}}\int_{\rnp}\g(y,z) F^1(z)dz$, for all $\phi\in L^{\rho}(\rn)$ we have that
	\begin{align*}
		\bigg\langle \mathcal{A}^q_q(VF),\varphi\bigg\rangle &= \int_{\rn}\int_0^{\infty}\int_{B_{y_n}(x')}\big|VF(y',y_n)\big|^q\dfrac{dy'dy_n}{y_n^{n}}\phi(x')dx'\\
		&\leq C\int_{\rn}\int_0^{\infty}\fint_{B_{y_n}(y')}\phi(x')dx'[y_n(G_1|F|)(y',y_n)]^qdy'dy_n\\
		&\leq C\int_{\rn}\int_0^{\infty}[y_n(G_1|F|)(y',y_n)]^qM_{y_n}\phi(y')dy'dy_n\\
		&\leq C\big\|(y',y_n)\mapsto y_nG_1|F|\big\|^q_{L^{p}L^q(\rnp)}\|M_{\cdot}\phi\|_{L^{\rho'}L^{\infty}(\rnp)}\\
		&\leq C\|F\|^{q}_{T^{r,\eta}_{-\frac{b+1/\eta}{n-1}}}\|\mathcal{M}\phi\|_{L^{\rho'}(\rn)}\\
		&\leq C\|F\|^{q}_{T^{\tau,\eta}_{s_{\eta}}}\|\phi\|_{L^{\rho'}(\rn)}.
	\end{align*}
	Note that the penultimate inequality follows from Remark \ref{rmk:Gb-weighted} with $r\in(\tau,p)$ such that $ \frac{1}{r}+\frac{1}{n-1}\leq\frac{1}{\tau}$ and $b=(n-1)(\frac{1}{\tau}-\frac{1}{r})$ while the last bound comes from \eqref{embedding-tents}. Collecting and summing up all the estimates on the $III_i$'s, we find that
	\begin{align*}
		\bigg\|\int_{\rnp}\g(\cdot,z) F(z)dz\bigg\|_{T^{p,q}_{\sigma_q}}\leq C\|F\|_{T^{\tau,\eta}_{s_{\eta}}}.
	\end{align*}
	The remaining estimate reads
	\begin{align*}
		\bigg\|\int_{\rnp}\nabla_z\g(\cdot,z) H(z)dz\bigg\|_{T^{p,q}_{\sigma_q}}\leq C\|H\|_{T^{\va,\sigma}_{s_{\sigma}}}.
	\end{align*}
	The argument we plan to use here is similar to the previous one. In fact, for $(y',y_n)\in \rnp$ we write
	\begin{align*}
		y_n\bigg|\int_{\rnp}\nabla_z\g(y,z) H(z)dz\bigg|\leq \sum_{k=1}^{3}\varGamma_k(y',y_n),
	\end{align*}
	with 
	\begin{align*} \varGamma_1(y',y_n)&=	y_n\int_{\rn\setminus B_{4y_n}(y')}\int_{0}^{\infty}|\nabla_{z}\g(y,z)| |H(z)|dz\\
		\varGamma_2(y',y_n)&=y_n\int_{B_{4y_n}(y')}\int_{4y_n}^{\infty}|\nabla_{z}\g(y,z)| |H(z)|dz\\
		\varGamma_3(y',y_n)&=y_n\int_{B_{4y_n}(y')}\int_{0}^{4y_n}|\nabla_{z}\g(y,z)| |H(z)|dz.
	\end{align*}
	Mimicking the proof of the Claim \ref{claim:1}, it is easy to see by Lemma \ref{lem:Green-pointwise-est} that for any $(x',y_n)\in \rnp$
	\begin{align*}
		\bigg(\fint_{B_{y_n}(x')} |\varGamma_1(y',y_n)|^qdy'\bigg)^{1/q}&\leq cG_1\bigg(\fint_{B_{y_n}(\cdot)} |H(z',\cdot)|dz'\bigg)(x',y_n).
	\end{align*}  
	Then, by applying Lemma \ref{lem:G-beta} with $\eta=\sigma$ and $\tau=\varLambda$, we deduce the  estimate
	\begin{equation*}
		\bigg\|\int_{\rnp}\nabla_z\g(\cdot,z) H^3(z)dz\bigg\|_{T^{p,q}_{\sigma_q}}\leq C \bigg\|G_1\bigg(\fint_{B_{y_n}(\cdot)} |H(z',\cdot)|dz'\bigg)\bigg\|_{L^{p}L^{q}(\rnp)}\leq C\|H\|_{T^{\varLambda,\sigma}_{s_{\sigma}}}.
	\end{equation*}
	Next, we show that
	\begin{align}\label{vargamma2}
		\bigg\|\int_{\rnp}\nabla_z\g(\cdot,z) H^2(z)dz\bigg\|_{T^{p,q}_{\sigma_q}}\leq C\|H\|_{T^{\va,\sigma}_{s_{\sigma}}}.	
	\end{align} 
	To achieve this, let us primarily observe that
	\begin{align*}
		|\varGamma_2(y',y_n)|&\leq Cy_n\big[\mathcal{A}^{4}_{\sigma}(z_n^{1/\sigma}H)(x')\big]\bigg(\int_{B_{4y_n}(x')}\int_{4y_n}^{\infty}\dfrac{z_n^{\frac{n-1}{\sigma-1}}dz_ndz'}{(|y'-z'|^2+|y_n-z_n|^2)^{\frac{n\sigma'}{2}}}\bigg)^{\frac{1}{\sigma'}}\\
		&\leq Cy_n^{1-\frac{1}{\sigma}}\big[\mathcal{A}^{4}_{\sigma}(z_n^{1/\sigma}H)(x')\big],\quad x'\in B(y',y_n).
	\end{align*}
	Taking the $\va$-power of both sides of the last inequality and integrating with respect to the variable $x'$ leads to
	\begin{align*}
		|\varGamma_2(y',y_n)|&\leq Cy_{n}^{1-\frac{1}{\sigma}-\frac{n-1}\va}\big\|\mathcal{A}^{4}_{\sigma}(z_n^{1/\sigma}H)\big\|_{L^{\va}(\rn)}\leq C\|H\|_{T^{\va,\sigma}_{s_{\sigma}}},\quad  y'\in B(x',y_n).
	\end{align*}
	Let $\delta>0$ to be determined later. The preceding inequalities imply for each $x'\in \rn$
	\begin{align*}
		\mathcal{A}_q\bigg(y_n^{1+1/q}\int_{\rnp}\nabla_z\g(\cdot,z) H^2(z)dz\bigg)(x')&=	\bigg(\int_{0}^{\infty}\fint_{B_{y_n}(x')}\big|\varGamma_2(y',y_n)\big|^qdy'dy_n\bigg)^{1/q}\\
		&\leq \bigg[\bigg(\int_{0}^{\delta}+\int_{\delta}^{\infty}\bigg)\fint_{B_{y_n}(x')}\big|\varGamma_2(y',y_n)\big|^qdy'dy_n\bigg]^{1/q}\\
		&\leq C\delta^{(n-1)(\frac{1}{\va}-\frac{1}{p})}[\mathcal{A}^4_{\sigma}(z_n^{1/\sigma}H)(x')]+\delta^{-\frac{n-1}{p}}\|H\|_{T^{\va,\sigma}_{s_{\sigma}}}. 
	\end{align*}
	Optimizing with respect to $\delta$, that is, choosing $\delta=\bigg(\|H\|_{T^{\va,\sigma}_{s_{\sigma}}}/\big[\mathcal{A}^{4}_{\sigma}(z_n^{1/\sigma}H)(x')\big]\bigg)^{\frac{\va}{n-1}}$ yields
	\begin{align*}
		\bigg(\int_{0}^{\infty}\fint_{B_{y_n}(x')}\big|\varGamma_2(y',y_n)\big|^qdy'dy_n\bigg)^{\frac{1}{q}}&\leq C\|H\|^{1-\frac{\va}{p}}_{T^{\va,\sigma}_{s_{\sigma}}}\big[\mathcal{A}^{4}_{\sigma}(z_n^{1/\sigma}H)(x')\big]^{\frac{\va}{p}},\quad x'\in \rn.
	\end{align*}
	Taking the $L^{p}$-norm on both sides and using Remark \ref{rmk:ind-aperture} gives \eqref{vargamma2}. Finally, one claims that the $T^{p,q}_{\sigma_q}$-norm of $\varGamma_3$ is controlled from above by a constant multiple of $\|H\|_{T^{\va,\sigma}_{s_{\sigma}}}$. This is derived from a simple duality argument in the same fashion as before. The proof of Proposition \ref{prop:Green-pot-bound} is now complete.
	
	We can now summarize the findings obtained above into a single theorem establishing the well-posedness of System \eqref{eq:S-eq} for boundary data in the scale of Triebel-Lizorkin space with negative amount of smoothness. 
	We say that a pair $(u,\pi)$ is a solution to \eqref{eq:S-eq} if $u$ and $\pi$ satisfy the relations 
	\begin{equation}\label{eq:int-eq}
		u(x)=\mathcal{H}f(x)+\mathscr{G}(F,H)(x),\quad
		\pi(x)=\mathcal{E}f(x)+\Psi(F,H)(x),\hspace{0.2cm}x\in \rnp.
	\end{equation}
	\begin{theorem}\label{thm3}
		Assume that the numbers $\eta,\tau,\sigma,\va$ and $p, q$ are as in Proposition \ref{prop:Green-pot-bound}. Then for any $f\in [\dot{F}^{-1/q}_{p,q}(\rn)]^{n}$, $F\in \X_{\tau,\eta}$ and $H\in \X_{\va,\sigma}$, the Stokes system \eqref{eq:S-eq} has a solution $(u,\pi)\in \X_{p,q}\times \Z_{p,q}$  (in the sense made precise in \eqref{eq:int-eq}) and the following estimate holds:
		\begin{equation}\|u\|_{\X_{p,q}}+\|\pi\|_{\Z_{p,q}}\leq C(\|f\|_{\dot{F}^{-\frac{1}{q}}_{p,q}(\rn)}+\|F\|_{\X_{\tau,\eta}}+\|H\|_{\X_{\va,\sigma}})
		\end{equation}	
		for some constant $C>0$ independent of $f$, $F$ and $H$.
	\end{theorem}

	\begin{remark}
		Practically, Theorem \ref{thm3} can easily be extended to the case where the vector field $u$ is not necessarily solenoidal, i.e. $\dv u=\phi$ under suitable conditions on $\phi$ $($integrability and compatibility$)$ using the  formulation derived in \cite[formula 2.32]{Solo}, see also \cite{Catta} so that our result gives an alternative approach to the Dirichlet problem for the Stokes system $($to be compared to \cite{ANR} and \cite{Fsa} wherein the analysis is carried out in weigthed Sobolev spaces and Lebesgue spaces, respectively$)$.
		These estimates of the velocity field and the pressure in weighted tent framework  against boundary data in low regularity spaces are new and generalize well-known results. In fact, our boundary class $\dot{F}^{-\frac{1}{q}}_{p,q}(\rn)$ contains the homogeneous Sobolev space $\dot{H}^{s,r}(\rn)$  for $-1/q<s<(n-1)/r$ under the conditions $1\leq r<p$ and $(n-1)/r-s=1/q+(n-1)/p$.  
	\end{remark}
	\section{Proofs of main results}\label{s:proofs}
	The proof of Theorem \ref{thm:unforced} can easily be deduced from that of Theorem \ref{thm1}. Hence, only the proofs of Theorems \ref{thm1} and \ref{thm:solvability-persistance} will be exposed in details and in the process, one essentially relies on preliminary results obtained in Section \ref{s:2}.
	
	\begin{proof}[Proof of Theorem \ref{thm1}]
		Let  $f\in [\dot{H}^{-\frac{1}{2},2(n-1)}(\rn)]^{n}$ with $n>2$ and assume $F\in \X_{\tau,\eta}$ for $\eta\in (1,2)$ and $\tau\in (\eta,2(n-1))$ with $\dfrac{1}{\eta}+\dfrac{n-1}{\tau}=3$. Equip the Banach space $\X\times \Z$ with the norm $\|\cdot\|:=\|\cdot\|_{\X}+\|\cdot\|_{\Z}$ and introduce the operators $\mathscr{L}$ defined by
		\[\mathscr{L}(u,\pi)=\big(\mathcal{H}f+\mathscr{G}[F,u\otimes u],\hspace{0.1cm} \mathcal{E}(f)+\Psi[F,u\otimes u]\big)\]
		where $\mathcal{H}$ and $\mathcal{E}$ are given by \eqref{eq:linear-int-eqs}. A solution of \eqref{eq:NS-eq}-\eqref{bc} according to the definition given by \eqref{eq:int-eq} is a couple $(u,\pi)$ satisfying the fixed point equation 
		\begin{align}\label{FPE}
			(u,\pi)=\mathscr{L}(u,\pi)\hspace{0.2cm}\mbox{in}\hspace{0.2cm}\rnp.   
		\end{align} 
		Using a Banach fixed point argument, we wish to show that the latter equation admits a solution in $\X\times \Z$. Assume that $(u,\pi)$ and $(v,\pi')$ in $\X\times \Z$ satisfy \eqref{FPE}  and use Proposition \ref{prop:Green-pot-bound} with $(p,q)=(2(n-1),2)$, $(\va,\sigma)=(n-1,1)$ to get	
		\begin{align}\label{eq:contraction}
			\nonumber\big\|\Lm(u,\pi)-\Lm(v,\pi')\big\|&=\big\|\mathscr{G}(0,u\otimes u-v\otimes v)\big\|_{\X}+\big\|\Psi(0,u\otimes u-v\otimes v)\big\|_{\Z}\\
			\nonumber& \leq C\|u\otimes u-v\otimes v\|_{\X_{n-1,1}}\\
			\nonumber&\leq C(\|u\otimes(u-v)\|_{\X_{n-1,1}}+\|(u-v)\otimes v\|_{\X_{n-1,1}})\\
			\nonumber&\leq C(\sup_{x_n>0}x_n^{2}\big\|[u\otimes(u-v)](\cdot,x_n)\big\|_{L^{\infty}(\rn)}+\|u\otimes (u-v)\|_{T^{n-1,1}_{s_1}}+\\
			\nonumber&\hspace{1.2cm}\sup_{x_n>0}x_n^{2}\|[(u-v)\otimes v](\cdot,x_n)\|_{L^{\infty}(\rn)}+\|(u-v)\otimes v\|_{T^{n-1,1}_{s_1}})\\
			\nonumber&\leq C\bigg(\sup_{x_n>0}x_n\|u(\cdot,x_n)\|_{L^{\infty}(\rn)}\sup_{x_n>0}x_n\|(u-v)(\cdot,x_n)\|_{L^{\infty}(\rn)}+\\
			\nonumber		&\hspace{1.3cm}\sup_{x_n>0}x_n\|(u-v)(\cdot,x_n)\|_{L^{\infty}(\rn)}\sup_{x_n>0}x_n\|v(\cdot,x_n)\|_{L^{\infty}(\rn)}+\\
			\nonumber&\hspace{1.82cm}\|u\|_{T^{2(n-1),2}_{s_2}}\|u-v\|_{T^{2(n-1),2}_{s_2}}+\|u-v\|_{T^{2(n-1),2}_{s_2}}\|v\|_{T^{2(n-1),2}_{s_2}}\bigg)\\
			&\leq C\|u-v\|_{\X}(\|u\|_{\X}+\|v\|_{\X}).
		\end{align}
		One can deduce from the previous estimate and in light of Lemma \ref{lem:lin-est} and Proposition \ref{prop:Green-pot-bound} the bound
		\begin{equation}\label{eq:self-map}
			\big\|\Lm(u,\pi)\big\|
			\leq C(\|u\|^{2}_{\X}+\|f\|_{\dot{H}^{-1/2,2(n-1)}(\rn)}+\|F\|_{\X_{\tau,\eta}}).
		\end{equation}
		Now pick $\varepsilon>0$ such that $\|f\|_{\dot{H}^{-1/2,2(n-1)}(\rn)} +\|F\|_{\X_{\tau,\eta}}\leq \varepsilon$. If $\varepsilon$ is sufficiently small, $\varepsilon<\min(1/4C,1/2)$ then it readily follows from \eqref{eq:contraction} and \eqref{eq:self-map} that $\Lm$ has a unique fixed point in a closed ball of $\E$ centered at the origin with radius $2\varepsilon$.
	\end{proof}	
	\begin{proof}[Proof of Theorem \ref{thm:solvability-persistance}]	
		Let $2< q<p<\infty$. Further, let $\eta_1\in (1,q)$ and $\tau_1\in (\eta_1,p)$ such that \begin{equation}\label{cds-exp}\frac{1}{\eta_1}+\frac{n-1}{\tau_1}=2+\frac{1}{q}+\frac{n-1}{p}.
		\end{equation}
		Assume $f\in \dot{H}^{-\frac{1}{2},2(n-1)}\cap \dot{F}^{-\frac{1}{q}}_{p,q}(\rn)$ and $F\in \X_{\tau_1,\eta_1}\cap \X_{\tau,\eta}$. 
		We remark that the solution found above may be realized as the unique limit in $\E$ of the following sequence of approximations given by
		\begin{equation*}
			\begin{cases}
				(u_{1},\pi_1)=(\mathcal{H}(f),\mathcal{E}(f))\\ 
				(u_{j+1},\pi_{j+1})=(\mathscr{G}[F,u_j\otimes u_j]+u_1,\Psi[F,u_j\otimes u_j]+\pi_1),\hspace{0.1cm} j=1,2,...  
			\end{cases}
		\end{equation*}
		Each element of this sequence belongs to $\X_{p,q}\times \Z_{p,q}$. In fact, since $(u_1,\pi_1)\in \X_{p,q}\times \Z_{p,q}$ (see Lemma \ref{lem:lin-est}) one may proceed via an induction argument to prove the claim. Choose $(\sigma,\va)$ such that $\frac{1}{\sigma}=\frac{1}{2}+\frac{1}{q}$, $\frac{1}{\va}=\frac{1}{2(n-1)}+\frac{1}{p}$ and invoke Proposition \ref{prop:Green-pot-bound}, H\"{o}lder's inequality in tent spaces simultaneously to have for each $j$, 
		\begin{align*}
			&\|(u_{j+1},\pi_{j+1})\|_{\X_{p,q}\times \Z_{p,q}}\\
			&=\|\mathscr{G}[F,u_{j}\otimes u_j]+u_1\|_{\X_{p,q}}+\|\Psi[F,u_j\otimes u_j]+\pi_j\|_{\Z_{p,q}}\\
			&\leq C\big(\|f\|_{\dot{F}^{-\frac{1}{q}}_{p,q}(\rn)}+\|F\|_{\X_{\tau_1,\eta_1}}+\|u_j\otimes u_j\|_{\X_{\va,\sigma}}\big)\\
			&\leq C\big(\|f\|_{\dot{F}^{-\frac{1}{q}}_{p,q}(\rn)}+\|F\|_{\X_{\tau_1,\eta_1}}+\|u_j\|_{T^{2(n-1),2}_{s_2}}\|u_j\|_{T^{p,q}_{s_q}}+\sup_{x_n>0}x^{\frac{1}{\sigma}+\frac{n-1}{\va}}_n\|u_j\otimes
			u_{j}(\cdot,x_n)\|_{L^{\infty}}\big)\\
			&\leq C\big(\|f\|_{\dot{F}^{-\frac{1}{q}}_{p,q}(\rn)}+\|F\|_{\X_{\tau_1,\eta_1}}+\|u_j\|_{\X}\|u_j\|_{\X_{p,q}}\big)
		\end{align*}
		so that if $(u_j,\pi_j)\in \X_{p,q}\times \Z_{p,q}$, then so is $(u_{j+1},\pi_{j+1})$. Next, we show that the latter sequence is Cauchy in $\X_{p,q}\times \Z_{p,q}$. We estimate $(w_j,q_j)=(u_{j+1}-u_j,\pi_{j+1}-\pi_j)$, $j=1,2,...$
		\begin{align*}
			\nonumber\|(w_{j},q_{j})\|_{\X_{p,q}\times \Z_{p,q}}&=\big\|\mathscr{G}[0,u_{j}\otimes u_{j}-u_{j-1}\otimes u_{j-1}]\big\|_{\X_{p,q}}+\big\|\Psi[0,u_{j}\otimes u_{j}-u_{j-1}\otimes u_{j-1}]\big\|_{\Z_{p,q}}\\
			\nonumber&\leq C\|u_{j}\otimes u_{j}-u_{j-1}\otimes u_{j-1}\|_{\X_{\va,\sigma}}\\
			\nonumber&\leq C\|u_{j}\otimes w_{j-1}+w_{j-1}\otimes u_{j-1}\|_{\X_{\va,\sigma}}\\
			&\leq C\|w_{j-1}\|_{\X_{p,q}}(\|u_j\|_{\X}+\|u_{j-1}\|_{\X})\\
			&\leq C\|(w_{j-1},q_{j-1})\|_{\X_{p,q}\times \Z_{p,q}}(\|u_j\|_{\X}+\|u_{j-1}\|_{\X}).
		\end{align*}
		Let $\varepsilon>0$ be as in Theorem \ref{thm1} and take $0<\varepsilon_{\star}<\varepsilon$. If $\|f\|_{\dot{H}^{-\frac{1}{2},2(n-1)}}+\|F\|_{\X_{\tau_1,\eta_1}}\leq \varepsilon_{\star}$, then the conclusion of Theorem \ref{thm1} shows that $\|u_j\|_{\X}\leq 2\varepsilon_{\star}$. Whence,
		\begin{align*}
			\|(w_{j},q_{j})\|_{\X_{p,q}\times \Z_{p,q}}\leq 4C\varepsilon_{\star}\|(w_{j-1},q_{j-1})\|_{\X_{p,q}\times \Z_{p,q}}.
		\end{align*}
		A simple iteration of the previous inequality yields
		\begin{align*}
			\|(w_{j},q_{j})\|_{\X_{p,q}\times \Z_{p,q}}\leq (4C\varepsilon_{\star})^{j-1}\|(w_{1},q_{1})\|_{\X_{p,q}\times \Z_{p,q}}
		\end{align*}
		thus implying the convergence of the sequence $(w_{j},q_{j})$ in $\X_{p,q}\times \Z_{p,q}$ since  $\varepsilon_{\star}<1/4C$. The limit of this sequence solves \eqref{FPE} and by uniqueness, it is the same as the solution constructed in Theorem \ref{thm1}.  
	\end{proof}

	\section{Appendix}
	Here we sketch the proof of Lemma \ref{lem:equiv-L^q}. Let $K\subset \rnp$ be a compact set. Then by Lemma \ref{lem:E_S} we know that $E(K)=\{x'\in \rn:K\cap \varGamma(x')\neq \emptyset\}$ has a finite Lebesgue measure. 
	\begin{proof}[Proof of Lemma \ref{lem:equiv-L^q}]
		Let us denote by $\mathbf{1}_K$ the characteristic function of the compact set $K$. 
		If $p\leq q$, then via H\"{o}lder's inequality, one obtains	
		\begin{align*}
			\|\mathbf{1}_Kf\|_{T^{p,q}_{s_q}}&=\bigg(\int_{\rn}\bigg(\iint_{\varGamma(x')}\mathbf{1}_K|f|^qy_n^{-(n-1)}dy'dy_n\bigg)^{p/q}dx'\bigg)^{1/p}\\
			&=\bigg(\int_{E(K)}\bigg(\iint_{\varGamma(x')}|f|^qy_n^{-(n-1)}dy'dy_n\bigg)^{p/q}dx'\bigg)^{1/p}\\
			&\leq \bigg(\int_{E(K)}\iint_{\varGamma(x')}|f|^qy_n^{-(n-1)}dy'dy_ndx'\bigg)^{1/q}\big|E(K)\big|^{\frac{1}{p}-\frac{1}{q}}\\
			&\leq \big|E(K)\big|^{\frac{1}{p}-\frac{1}{q}}\bigg(\int_{\rn}\iint_{\varGamma(x')\cap K}|f|^qy_n^{-(n-1)}dy'dy_ndx'\bigg)^{1/q}\\
			&\leq C\big|E(K)\big|^{\frac{1}{p}-\frac{1}{q}}\|f\|_{L^{q}(K)}.
		\end{align*}	
		Moving on, for $q<p$, applying Minkowski's inequality implies  
		\begin{align*}
			\|\mathbf{1}_Kf\|_{T^{p,q}_{s_q}}&=\bigg(\int_{\rn}\bigg(\iint_{\varGamma(x')}\mathbf{1}_K|f|^qy_n^{1-n}dy'dy_n\bigg)^{p/q}dx'\bigg)^{1/p}\\
			&=\bigg(\int_{\rn}\bigg(\int_{\rn}\int^{\infty}_0\mathbf{1}_{B_{y_n}(y')}(x')\mathbf{1}_K(y',y_n)|f|^qy_n^{1-n}dy_ndy'\bigg)^{p/q}dx'\bigg)^{1/p}\\
			&\leq C_K\bigg(\int_{\rn}\int_{0}^{\infty}\mathbf{1}_K|f|^qdy'dy_n\bigg)^{1/q}\\
			&\leq C_K\|f\|_{L^{q}(K)}.
		\end{align*}	
		Assuming that $p\leq q$, we use Lemma \ref{lem:Aver} and Minkowski's inequality simultaneously to get	
		\begin{align*}
			\|f\|_{L^q(K)}&=\bigg(\int_{\rnp}\mathbf{1}_K|f|^qdy'dy_n\bigg)^{1/q}\\
			&\leq C\bigg(\int_{\rn}\iint_{\varGamma(x')}\mathbf{1}_K|f|^qy_n^{1-n}dy'dy_ndx'\bigg)^{1/q}\\
			&\leq C_K\bigg(\int_{\rn}\int_0^{\infty}\bigg(\int_{E(K)}\mathbf{1}_{B_{y_n}(y')}(x')y_n^{\frac{(n-1)p}{q}}|f|^pdx'\bigg)^{q/p}dy'dy_n\bigg)^{1/q}\\
			&\leq C_K\bigg(\int_{\rn}\bigg(\iint_{\varGamma(x')}|f|^qy_n^{1-n}dy'dy_n\bigg)^{p/q}dx'\bigg)^{1/p}\\
			&\leq C_K\|f\|_{T^{p,q}_{s_q}}.
		\end{align*}	
		When $p>q$, the desired bound follows from  H\"{o}lder's inequality. Indeed, we have	
		\begin{align*}
			\|f\|_{L^q(K)}&=\bigg(\int_{\rnp}\mathbf{1}_K|f|^qdy'dy_n\bigg)^{1/q}\\
			&\leq C\bigg(\int_{\rn}\iint_{\varGamma(x')}\mathbf{1}_K|f|^qy_n^{1-n}dy'dy_ndx'\bigg)^{1/q}\\
			&\leq C\bigg(\int_{E(K)}\iint_{\varGamma(x')}|f|^qy_n^{1-n}dy'dy_ndx'\bigg)^{1/q}\\
			&\leq C\bigg(\int_{\rn}\bigg(\iint_{\varGamma(x')}|f|^qy_n^{1-n}dy'dy_n\bigg)^{p/q}dx'\bigg)^{1/p}\big|E(K)\big|^{\frac{1}{q}-\frac{1}{p}}\\
			&\leq C\big|E(K)\big|^{\frac{1}{q}-\frac{1}{p}}\|f\|_{T^{p,q}_{s_q}}.
		\end{align*}	
	\end{proof} 
	Next, we prove Lemma \ref{lem:Fourierkernels}. Before we proceed, let us recall some properties of the Fourier transform. Let $\varphi\in \mathcal{S}(\rn)$, its Fourier transform is denoted by the standard notation  $\widehat{\varphi}$ (or $\mathcal{F}\varphi$) and defined as $\widehat{\varphi}(\xi')=\int_{\rn}e^{-iy'\cdot \xi'}\varphi(y')dy'$ for all $\xi'\in \rn$. If $\varphi^{\vee}=\varphi(-\cdot)$ and  $u$ is a tempered distribution in $\rn$, then the map $u^{\vee}$ defined by $u^{\vee}:\mathcal{S}(\rn)\rightarrow \R$, $u^{\vee}(\varphi)=\big<u,\varphi^{\vee}\big>$ for all $\varphi\in  \mathcal{S}(\rn)$ is a tempered distribution in $\rn$ and 
	\begin{equation}\label{eq:Fourier-flip}
		\mathcal{F}^{-1}u=(2\pi)^{1-n}\widehat{u^{\vee}}.
	\end{equation} 
	In addition to this formula, we will systematically use below the  Fourier transforms of the integrable functions $f(z)=(z^2+b^2)^{-1}$, $z\in \R$, $b>0$ and $p(y')=\frac{2}{\omega_{n-1}}(1+|y'|^2)^{-n/2}$, $y'\in \rn$  given respectively by 
	\begin{equation}\label{eq:fourier-specialfunct}
		\widehat{f}(\eta)=\frac{\pi}{b}e^{-b|\eta|},\,\,\, \forall\, \eta\in \R \,\,\mbox{and}\,\, \mbox{in}\,\, \mathcal{S}'(\R)\,\,\,\mbox{and}\,\,\,\,  \widehat{p}(\xi')=e^{-|\xi'|},\,\, \xi'\in \rn\setminus\{0\}.  
	\end{equation}
	\begin{proof}[Proof of Lemma \ref{lem:Fourierkernels}]
		We start with the computation of $\widehat{P_0}$ by means of partial Fourier transform. Observe that $ \mathcal{F}_{x'}(K_{nn}(\cdot,x_n))(\xi')=-x_n|\xi'|e^{-x_n|\xi'|}-e^{-x_n|\xi'|}$ $\forall \,\xi'\in \rn\setminus\{0\}$.
		This is easily seen from the formula $K_{nn}(x)=\frac{2}{\omega_{n-1}}\big(x_n\partial_n\frac{x_n}{|x|^n}-\frac{x_n}{|x|^n}\big)$ and \eqref{eq:fourier-specialfunct}. From the identity $K_{nn}(x)=x_{n}^{-(n-1)}P_0(\frac{x'}{x_n})$, $x=(x',x_n)\in \rnn_+$, one deduces that for each $\xi'\in \rn$
		\begin{equation}
			\widehat{P_{0}}(\xi')=-|\xi'|e^{-|\xi'|}-e^{-|\xi'|}.    
		\end{equation}
		The Fourier transform of $Q_{jk}$, $j,k=1,...,n-1$ can be derived from the simple identity \begin{equation}\label{eq:Qjk}Q_{jk}(x')=-\delta_{jk}p(x')-\frac{2}{\omega_{n-1}(n-2)}\partial^2_{jk}(|x'|^2+1)^{\frac{2-n}{2}}
		\end{equation}
		if $\widehat{(|\cdot|^2+1)^{\frac{2-n}{2}}}$ is explicit. In this regard, note that $|x|^{-(n-2)}$, $x\in \rnn\setminus\{0\}$ belongs to $\mathcal{S}'(\rnn)$ and for each $\xi'\in \rn\setminus\{0\}$ and $\eta\in \R$ we have \[\widehat{|\cdot|^{-(n-2)}}(\xi',\eta)=(n-2)\omega_{n-1}(|\xi'|^2+|\eta|^2)^{-1}. \] 
		Hence, using \eqref{eq:Fourier-flip} in one dimension and \eqref{eq:fourier-specialfunct} with $b=|\xi'|$ we arrive at
		\begin{align*}
			\mathcal{F}_{x'}(|(x',x_n)|^{2-n})(\xi')&=\mathcal{F}^{-1}_{\eta}[\widehat{|\cdot|^{-(n-2)}}(\xi',\eta)](x_n)\\
			&=(n-2)\omega_{n-1}\mathcal{F}^{-1}_{\eta}\bigg(\frac{1}{|\xi'|^2+\eta^2}\bigg)(x_n)\\
			&=(2\pi)^{-1}(n-2)\omega_{n-1}\mathcal{F}_{\eta}\bigg(\frac{1}{|\xi'|^2+\eta^2}\bigg)(x_n)=\frac{(n-2)\omega_{n-1}}{2}|\xi'|^{-1}e^{-|x_n\xi'|}.
		\end{align*}
		In particular, when $x_n>0$, $\mathcal{F}_{x'}(|x|^{2-n})(\xi')=\frac{(n-2)\omega_{n-1}}{2}|\xi'|^{-1}e^{-x_n|\xi'|}$ for all $\xi'\in\rn\setminus\{0\}$. As a consequence, $\widehat{(|\cdot|^2+1)^{\frac{2-n}{2}}}=\frac{(n-2)\omega_{n-1}}{2}|\xi'|^{-1}e^{-|\xi'|}$, $\xi'\in \rn\setminus\{0\}$ and by \eqref{eq:fourier-specialfunct}, it follows from \eqref{eq:Qjk} that
		\begin{equation*}\widehat{Q_{jk}}(\xi')=-\delta_{jk}e^{-|\xi'|}+\frac{\xi'_{j}\xi'_k}{|\xi'|}e^{-|\xi'|},\,\,\,\xi'\in \rn\setminus\{0\}.
		\end{equation*}
		Finally, since $P_j(x')=\partial_jp(x')$, one uses \eqref{eq:fourier-specialfunct} to get $\widehat{P_j}(\xi)=i\xi'_je^{-|\xi'|}$, $\xi'\in \rn$.
	\end{proof}
	We now turn to the proof of Proposition \ref{prop:Stokes-Triebel}. Given $\lambda>0$, a kernel $\psi$ defined on $\rn$ and a distribution $f\in \mathcal{S}'(\rn)$ such that $\psi_{t}f=\psi_t\ast f$, $t>0$ is continuous, the Peetre maximal function associated to $f$ is defined as
	\begin{equation}\label{Peetre}
		(\psi_tf)^{\ds}(y')=\sup_{z'\in \rn}|\psi_tf(z')|(1+t^{-1}|z'-y'|)^{-\lambda},\,\,y'\in \rn.   
	\end{equation}
	\begin{proof}[Proof of Proposition \ref{prop:Stokes-Triebel}]
		Let $f\in \dot{F}^{s}_{p,q}(\rn)$ with $s<0$  and $\psi\in \mathcal{A}_{\varLambda,m,r}$ where $r>s$, $m+s>\lambda$ with $\lambda\in (\va,[\va]+1)$, $\va=\max\{(n-1)/p,(n-1)/q\}$. Then $\psi_{x_n}f$ is a bounded continuous function. We primarily estimate the conical functional of the underlying convolution by the Littlewood-Paley functional of the Peetre maximal function associated to $f$.  Indeed, 
		\begin{align*}
			\mathcal{A}_{q}(y_n^{-s}\psi_{y_n}f)(x')&=\bigg(\iint_{\varGamma(x')} \big|y_n^{-s}\psi_{y_n}f(y')\big|^qy^{-n}_ndy'dy_n\bigg)^{1/q}\\   
			&\leq 2^{\lambda/q}\bigg(\int_{0}^{\infty}\int_{B_{y_n}(x')}|y_n^{-s}(\psi_{y_n}f)^{\ds}(x')|^qy_n^{-n}dy'dy_n\bigg)^{1/q}\\
			&\leq 2^{\lambda/q}\bigg(\int_{0}^{\infty}|y_n^{-s}(\psi_{y_n}f)^{\ds}(x')|^q\frac{dy_n}{y_n}\bigg)^{1/q}
		\end{align*}
		so that
		\begin{equation}\label{eq:tenttoPee}
			\|\psi_{y_n}f\|_{T^{p,q}_{s/(n-1)}}    \leq 2^{\lambda/q}\bigg\|\bigg(\int_0^{\infty}\big|y_n^{-s}(\psi_{y_n}f)^{\ds}(\cdot)\big|^q\frac{dy_n}{y_n}\bigg)^{1/q}\bigg\|_{L^p(\rn)}.
		\end{equation}
		Now, our goal is to estimate the function \[\chi(y,z')=|(\psi_{y_n}f)(z')|(1+y_n^{-1}|z'-y'|)^{-\lambda},\quad y\in \rnp,\,z'\in \rn.\] 
		After possibly subtracting a suitable polynomial to $f$, one may use the Calder\'{o}n reproducing formula for elements in $\dot{F}^{s}_{p,q}(\rn)$ to obtain the pointwise identity (see \cite[Theorem 3.1]{BC})
		\[\psi_{y_n}f(z')=\sum_{k\in \mathbb{Z}}\phi_k\ast \phi_k\ast \psi_{y_n}\ast f(z'),\,\,z'\in \rn\]
		where $\phi_k$ is as in Remark \ref{rmk:disccharac}. Let $l\in \mathbb{Z}$ and assume that $y_n\in [2^{-l-1},2^{-l}]$. Then, we have
		\begin{align*}
			\chi(y,z') &\leq \sum_{k\in \mathbb{Z}}|\psi_{y_n}\ast \phi_k\ast \phi_k\ast f(z')|(1+y_n^{-1}|y'-z'|)^{-\lambda}\\
			&\leq  \sum_{k\in \mathbb{Z}}\int_{\rn}|\psi_{y_n}\ast \phi_k(v')||\phi_k\ast f(z'-v')|dv'(1+y_n^{-1}|y'-z'|)^{-\lambda}\\
			&\leq 2^{(n-1)l} \sum_{k\in \mathbb{Z}}\phi^{\star}_kf(y')\int_{\rn}\big|\psi_{2^ly_n}\ast \phi_{k-l}(2^lv')\big|\frac{(1+2^{k}|y'-z'+v'|)^{\lambda}}{(1+2^{l}|y'-z'|)^{\lambda}}dv'\\
			&\leq \sum_{k\in \mathbb{Z}}\phi^{\star}_kf(y')\int_{\rn}\big|\psi_{2^ly_n}\ast \phi_{k-l}(v')\big|\frac{(1+2^{k-l}|2^l(y'-z')+v'|)^{\lambda}}{(1+2^{l}|y'-z'|)^{\lambda}}dv'.
		\end{align*}
		This clearly implies the pointwise estimate 
		\begin{equation}\label{peetrepointest}
			y_n^{-s}(\psi_{y_n}f)^{\ds}(y')    \leq \sum_{k\in \mathbb{Z}}a_{l-k}2^{sk}\phi^{\star}_kf(y');\quad y'\in \rn,\,\, y_n\in [2^{-l-1},2^{-l}]
		\end{equation}
		where the sequence $(a_k)_{k\in \mathbb{Z}}$ reads 
		\begin{equation*}
			a_k=2^{ks}\sup_{2^{-1}\leq \varrho\leq 1}\sup_{z'\in \rn}\int_{\rn}\big|\psi_{\varrho}\ast\phi_{-k}(v')\big|\dfrac{(1+2^{-k}|z'+v'|)^{\lambda}}{(1+|z'|)^{\lambda}}dv'.   
		\end{equation*}
		Hence, for any $y'\in \rn$, it follows from \eqref{peetrepointest} that
		\begin{align*}
			\bigg(\int_0^{\infty}\big[y_n^{-s}(\psi_{y_n} f)^{\ds}(y')\big]^qy_n^{-1}dy_n\bigg)^{1/q}& \leq  \bigg(\sum_{l\in \mathbb{Z}}\int_{2^{-l-1}}^{2^{-l}}\big[y_n^{-s}(\psi_{y_n}f)^{\ds}(y')\big]^qy_n^{-1}dy_n\bigg)^{1/q} \\
			&\leq \bigg(\sum_{l\in \mathbb{Z}}\bigg(\sum_{k\in \mathbb{Z}}a_{l-k}2^{ks}\phi^{\star}_kf(y')\bigg)^q\bigg)^{1/q}.
		\end{align*}
		The desired estimate will now follow from the above estimate, \eqref{eq:tenttoPee} and Remark \ref{rmk:disccharac} as soon as $(a_k)_k\in \ell^1$. Observe that 
		\begin{equation*}
			\dfrac{1+2^{-k}|z'+v'|}{1+|z'|}\leq \begin{cases}
				2^{-k}(1+|v'|)\,\,\,\,\mbox{if} \,\,\,\, k<0\\
				(1+2^{-k}|v'|)\,\,\,\, \mbox{if} \,\,\,\, k\geq 0
			\end{cases}   
		\end{equation*}
		and using the fact that $\psi$ belongs to $\mathcal{A}_{\va,m,r}$,  an application of Lemma 2.4 in \cite{BC} yields the estimate 
		\begin{equation}
			|\psi\ast \phi_{-k}(v')|\leq \begin{cases}
				2^{km}(1+|v'|)^{-(n+[\varLambda])}\,\,\,\,\,\,\,\,\,\,\,\,\,\,\,\,\,\,\,\,\,\,\,\,\,\,\,\,\,\,\, \mbox{if} \,\,\,\, k< 0\\
				2^{-k(n-1+r)}(1+2^{-k}|v'|)^{-(n+[\varLambda])}\,\,\,\,\mbox{if} \,\,\,\, k\geq 0.
			\end{cases}  
		\end{equation}
		Thus, by definition of $a_k$ and after rescaling, one obtains for each $k\in \mathbb{Z}$ the estimate 
		\begin{equation*}
			a_k\leq C\begin{cases}
				2^{-k(\lambda-m-s)} \,\,\,\, \mbox{if} \,\,\,\, k< 0\\
				2^{-k(r-s)}  \,\,\,\,\,\,\,\,\,\,\,\, \mbox{if} \,\,\,\, k\geq 0
			\end{cases}    
		\end{equation*}
		from which it easily follows by the assumptions on the parameters $\va,m,r$ that $(a_k)_k\in \ell^1$. The proof of Proposition \ref{prop:Stokes-Triebel} is now complete.
	\end{proof}
	
	\bibliographystyle{acm}

\end{document}